\documentclass[12pt]{article}
\usepackage{amsmath,amssymb}
\usepackage{graphicx}
\begin{document}

\title{Dynamic Stability of the 3D Axi-symmetric Navier-Stokes Equations
with Swirl}
\author{Thomas Y. Hou\thanks{Applied and Comput. Math, Caltech, Pasadena,
CA 91125. Email: hou@acm.caltech.edu.} \and Congming Li \thanks{Department of
Applied Mathematics, University of Colorado, Boulder, CO. 80309. Email: cli@colorado.edu}}

\maketitle

\begin{abstract}
In this paper, we study the dynamic stability of the
3D axisymmetric Navier-Stokes Equations with swirl. To this purpose, 
we propose a new one-dimensional (1D) model which approximates the
Navier-Stokes equations along the symmetry axis. An important
property of this 1D model is that one can construct from its solutions a family 
of exact solutions of the 3D Navier-Stokes equations. The nonlinear structure of the 1D
model has some very interesting properties. On one hand, it can lead to tremendous
dynamic growth of the solution within a short time. On the 
other hand, it has a surprising dynamic 
depletion mechanism that prevents the solution from blowing up 
in finite time. By exploiting this special nonlinear structure,
we prove the global regularity of the 3D Navier-Stokes equations
for a family of initial data, whose solutions can lead to
large dynamic growth, but yet have global smooth solutions.
\end{abstract}

\section{Introduction.}

Despite a great deal of effort by many mathematicians and physicists,
the question of whether the solution of the 3D Navier-Stokes equations can develop 
a finite time singularity from a smooth initial condition with finite
energy remains one of the most outstanding open problems \cite{Clay}. 
A main difficulty in obtaining the global regularity of the
3D Navier-Stokes equations is due to the presence of the vortex
stretching, which is absent for the 2D problem. Under suitable
smallness assumption on the initial condition, global
existence and regularity results have been obtained for some
time \cite{Lady70,CF88,Temam01,BM02}. But these methods based on energy
estimates do not generalize to the 3D Navier-Stokes with large 
data. Energy estimates seem to be too crude to give a definite answer to 
whether diffusion is strong enough to control the nonlinear 
growth due to vortex stretching. A more refined analysis 
which takes into account the special nature of the nonlinearities
and their local interactions seems to be needed.

In this paper, we study the dynamic stability property of the
3D axisymmetric Navier-Stokes Equations with swirl. We show that
there is a very subtle dynamic depletion mechanism of vortex stretching
in the 3D Navier-Stokes equations. On one hand, the nonlinear
vortex stretching term is responsible for producing a large
dynamic growth in vorticity in early times. On the other hand, 
the special structure of the nonlinearity can also lead to 
dynamic depletion and cancellation of vortex stretching, thus 
avoiding the finite time blowup of the Navier-Stokes equations. 

This subtle nonlinear stability property can be best illustrated
by a new 1D model which we introduce in this paper. 
This 1D model approximates the 3D axisymmetric Navier-Stokes equations
along the symmetry axis. By the well-known Caffarelli-Kohn-Nirenberg 
theory \cite{CKN82} (see also \cite{Lin98}), 
the singularity set of any suitable weak solution of the 3D 
Navier-Stokes equations has one-dimensional Hausdorff measure zero.
In the 
case of axisymmetric 3D Navier-Stokes equations with swirl, if 
there is any singularity, it must be along the symmetry axis. 
Thus it makes sense to focus our effort to understand the possible 
singular behavior of the 3D Navier-Stokes equations near the 
symmetry axis at $r=0$. By expanding the angular 
velocity ($u^\theta$), the angular vorticity ($\omega^\theta$),
and the angular stream function ($\psi^\theta$) around $r=0$,
we obtain the following coupled nonlinear partial differential
equations (see Section 2 for detailed derivations):
\begin{eqnarray}
&&\left( u_1 \right )_t + 2 \psi_1
\left ( u_1 \right )_z  = \nu ( u_1)_{zz} +
2 \left (\psi_1 \right )_z u_1 \label{eqn-u100} \\
&&\left( \omega_1 \right )_t + 2 \psi_1 \left ( \omega_1 \right )_z  =
\nu ( \omega_1)_{zz} + \left (u_1^2 \right )_z \label{eqn-w100}\\
&& - (\psi_1)_{zz} = \omega_1, \label{eqn-psi100}
\end{eqnarray}
where $u_1(z,t) \approx ( u^\theta)_r |_{r=0}$,
$\omega_1(z,t) \approx ( \omega^\theta )_r |_{r=0}$, and
$\psi_1(z,t) \approx ( \psi^\theta)_r |_{r=0}$.

What we find most surprising is that one can construct a family of
{\it exact solutions} from the above 1D model. Specifically, if 
$(u_1, \omega_1, \psi_1)$ is a solution of the 1D model 
(\ref{eqn-u100})-(\ref{eqn-psi100}), then 
\[
u^\theta (r,z,t) = r u_1(z,t), \quad 
\omega^\theta (r,z,t) = r \omega_1(z,t), \quad 
\psi^\theta (r,z,t) = r \psi_1(z,t), 
\]
is an exact solution of the 3D axisymmetric Navier-Stokes equations.
Thus the 1D model captures some essential nonlinear features
of the 3D Navier-Stokes equations. 
Further, if we let $\tilde u = u_1$, $\tilde v = -(\psi_1)_z$, and
$\tilde \psi = \psi_1$, then the 1D model can be rewritten as
\begin{eqnarray}
&& (\tilde{u})_t + 2 \tilde{\psi} (\tilde{u})_z  =
\nu ( \tilde{u})_{zz}
- 2 \tilde{v} \tilde{u}
\label{eqn-u20} \\
&& (\tilde{v})_t + 2 \tilde{\psi} (\tilde{v})_z  =
\nu (\tilde{v})_{zz} + (\tilde{u})^2 - (\tilde{v})^2
+c(t),
\label{eqn-v20}
\end{eqnarray}
where $\tilde {\psi}_z= -\tilde{v}$ and $c(t)$ is an integration
constant to ensure that  $\int \tilde{v} dz = 0$.
We will show that if the initial value of $\tilde u$ is small, but
$\tilde v$ is large and negative, then the solution of the 1D model can 
experience large growth. On the other hand, we also 
find a surprising dynamic depletion mechanism of nonlinearities
that prevents the 
solution from blowing up in finite time. This subtle nonlinear
cancellation is partly due to the special nature of the nonlinearities,
i.e. $-2 \tilde u \tilde v$ in (\ref{eqn-u20}), and 
$\tilde u^2 - \tilde v^2$ in (\ref{eqn-v20}). If one modifies 
the sign of the nonlinear term from $-2 \tilde u \tilde v$
to $2 \tilde u \tilde v$ or changes $\tilde u^2 - \tilde v^2$ to 
$\tilde u^2 + \tilde v^2$ or even modifies the
coefficient from $-2 \tilde u \tilde v$  to $-0.9 \tilde u \tilde v$,
the dynamic depletion mechanism can be changed completely.
Another interesting fact is that the convection term also helps to stabilize
the solution. It cancels some of the destabilizing terms from the
right hand side when we estimate the solution in a high order norm.
Specifically, we find that there is a miraculous cancellation
of nonlinear terms in the equation that governs 
the nonlinear quantity, $\tilde u_z^2 + \tilde v_z^2 $, i.e.
\begin{equation}
\left ({\tilde u}_z^2+{\tilde v}_z^2 \right )_t +
2 \tilde \psi \left (\tilde u_z^2 + {\tilde v}_z^2 \right )_z
 = \nu \left ( \tilde u_z^2  + \tilde v_z^2 \right )_{zz}
- 2 \nu \left [ (\tilde u_{zz})^2 + (\tilde v_{zz})^2  \right ] .
\end{equation}
Therefore, $\tilde u_z^2 + \tilde v_z^2 $ has a maximum principle.
This pointwise {\it a priori} estimate plays an essential role
in obtaining global regularity of the 1D model with or without 
viscosity. If one attempts to prove global regularity of the 
1D model using energy estimates, one cannot take full advantage
of this local cancellation of nonlinearities and would run into 
similar difficulties that we encounter for the 3D Navier-Stokes equations.

Finally, we construct a family of globally smooth solutions of the
3D Navier-Stokes equations with large initial data of finite
energy by using the solution of the 1D model. Specifically, we
look for the solution of the form:
\begin{eqnarray}
\tilde{u}^\theta &=& r \left (\bar{u}_1 (z,t) \phi (r) + u_1 (r,z,t) \right )\\
\tilde{\omega}^\theta &= &r \left (\bar{\omega}_1 (z,t) \phi (r) +
\omega_1 (r,z,t) \right ) \\
\tilde{\psi}^\theta &= &r \left (\bar{\psi}_1 (z,t) \phi (r) + \psi_1 (r,z,t)
\right ),
\end{eqnarray}
where $\bar{u}_1 $, $\bar{\omega}_1$ and $\bar{\psi}_1$ are 
solutions of the 1D model, $\phi (r)$ is a cut-off function to
ensure that the solution has finite energy. By using the 
{\it a priori} estimate of the solution of the 1D model and using
a delicate analysis, we prove that there exists a family of 
globally smooth functions $u_1(r,z,t)$, $\omega_1(r,z,t)$ and
$\psi_1 (r,z,t)$, such that $\tilde{u}^\theta$, 
$\tilde{\omega}^\theta$ and $\tilde{\psi}^\theta$ are solutions
of the 3D axisymmetric Navier-Stokes equations. Unlike the
other known global solutions with small data, the solutions
that we construct above using the 1D model can have large dynamic
growth for early times, which is induced by the dynamic growth of 
the corresponding solution of the 1D model, but yet 
the solution remains smooth for all times.

There has been some interesting development in the study of the
3D incompressible Navier-Stokes equations and related models. 
In particular, by exploiting the special structure of the governing
equations, Cao and Titi~\cite{CT05} prove the global well-posedness 
of the 3D viscous primitive equations which model large scale ocean 
and atmosphere dynamics. By taking advantage of the limiting
property of some rapidly oscillating operators and using 
non-linear averaging, Babin, Mahalov and Nicolaenko \cite{BMN01} 
prove existence on infinite time intervals 
of regular solutions to the 3D Navier-Stokes equations for 
some initial data characterized by uniformly large vorticity. 
Some interesting
progress has been made on the regularity of the axisymmetric
solutions of the Navier-Stokes equations, see e.g. \cite{CL02} 
and the references cited there. The 2D Boussinesq equations are
closely related to the 3D axisymmetric Navier-Stokes equations
with swirl (away from the symmetry axis). 
Recently, Chae \cite{Chae05} and Hou-Li \cite{HL05}
have proved independently the global existence of the 2D viscous
Boussinesq equations with viscosity entering only in the
fluid equation, but the density equation remains inviscid.
Recent studies by Constantin-Fefferman-Majda \cite{CFM96}
and Deng-Hou-Yu \cite{DHY05,DHY06} show that the local geometric 
regularity of the unit vorticity vector can play an important 
role in depleting vortex stretching dynamically. Motivated by 
these theoretical results, Hou and R. Li \cite{HL06b} have
recently re-investigated the well-known computations by Kerr 
\cite{Kerr93} for two anti-parallel vortex tubes, in which a 
finite time singularity of the 3D incompressible Euler equations 
was reported. The results of Hou and Li show that there is 
tremendous dynamic cancellation in the vortex stretching 
term due to local geometric regularity of the vortex lines.
Moreover, they show that the vorticity does not grow faster 
than double exponential in time and the velocity field remains 
bounded up to $T=19$, beyond the singularity time alleged in \cite{Kerr93}. 
Finally, we would like to mention the recent work of Gibbon et al
(see \cite{GHKR06} and the references therein) where they reveal
some interesting geometric properties of the Euler equations in 
quaternion-frames. 

The rest of the paper is organized as follows. In Section 2, we will
derive the 1D model for the 3D axisymmetric Navier-Stokes equations.
We discuss some of the properties of the 1D model in Section 3 and
prove the global existence of the inviscid 1D model using the Lagrangian 
coordinate. Section 4 is devoted to prove the global regularity of the full
1D model in the Eulerian coordinate. Finally in Section 5, we 
use the solutions of the 1D model to construct a family of solutions
of the 3D Navier-Stokes equations and prove that they remain smooth
for all times.

\section{Derivation of the 1D Model}

Consider
the 3D axi-symmetric incompressible Navier-Stokes equations with swirl.
\begin{equation}
  \left\{ \begin{array}{l}
    u_t + (u \cdot \nabla ) u = - \nabla p + \nu \Delta u,\\
    \nabla \cdot u = 0,\\
    u|_{t = 0} = u_0 ( \vec{x}), \quad \vec{x} = (x,y,z).
  \end{array} \label{nse} \right.
\end{equation}

Let 
\[
e_r = \left ( \frac{x}{r}, \frac{y}{r}, 0 \right ),
\;\;
e_\theta = \left (-\frac{y}{r}, \frac{x}{r}, 0 \right ),\;\;
e_z = \left ( 0, 0, 1 \right ), 
\]
be three unit vectors along the radial,
the angular, and the $z$ directions respectively,
$r=\sqrt{x^2+y^2}$. We will decompose the
velocity field as follows:
\begin{equation}
{\vec u} = v^r(r,z,t) e_r + u^\theta (r,z,t) e_\theta + v^z (r,z,t) e_z.
\end{equation}
In the above expression, $u^\theta$ is called the {\it swirl} component
of the velocity field ${\vec u}$.
The vorticity field can be expressed similarly
\begin{equation}
{\vec \omega } = -(u^\theta )_z (r,z,t) e_r + \omega^\theta (r,z,t) e_\theta +
\frac{1}{r} (r u^\theta )_r (r,z,t) e_z,
\end{equation}
where $\omega^\theta = v_z^r - v^z_r $.

To simplify our notation, we will use $u$ and $\omega$ to denote the 
angular velocity and vorticity 
components respectively, dropping the $\theta$ superscript in the rest
of the paper. One can 
derive evolution equations for $u$ and $\omega$ as follows 
(see e.g. \cite{BM02,CL02}).
\begin{eqnarray}
&& u_t + v^r u_r + v^z u_z = \nu \left ( \nabla^2 - \frac{1}{r^2} \right )
u  - \frac{1}{r} v^r u ,\label{eqn-u} \\
&& \omega_t + v^r \omega_r + v^z \omega_z
 = \nu \left ( \nabla^2 - \frac{1}{r^2} \right )
\omega +\frac{1}{r} \left ( u^2 \right )_z + \frac{1}{r}
v^r \omega ,\label{eqn-w} \\
&& - \left ( \nabla^2 - \frac{1}{r^2} \right ) \psi = \omega ,
\label{eqn-psi}
\end{eqnarray}
where $\psi$ is the angular component of the stream function, 
$v^r$ and $v^z$ can be expressed in terms of the angular stream function
$\psi$ as follows:
\begin{equation}
v^r = - \frac{\partial \psi}{\partial z}, \quad
v^z =  \frac{1}{r} \frac{\partial }{\partial r} (r \psi ),
\label{equ-v3r}
\end{equation}
and
$\nabla^2$ is defined as
\begin{equation}
\nabla^2 = \partial_r^2 + \frac{1}{r} \partial_r + \partial_z^2 .
\label{laplacian}
\end{equation}
Note that equations (\ref{eqn-u})-(\ref{eqn-psi}) completely
determine the evolution of the 3D axisymmetric Navier-Stokes
equations once the initial condition is given. 

Now, we will derive the 1D model for the 3D axisymmetric Navier-Stokes
equations.  By the well-known Caffarelli-Kohn-Nirenberg theory \cite{CKN82},
the singularity set of any suitable weak solution of the 3D Navier-Stokes 
equations has one-dimensional Hausdorff measure zero.
Thus, in the case of axisymmetric 3D Navier-Stokes
equations with swirl, if there is any singularity, it must be along the
symmetry axis, i.e. the $z$-axis. Therefore, we should focus our effort
to understand the possible singular behavior of the 3D Navier-Stokes
equations near the symmetry axis at $r=0$.

As observed by Liu and Wang in \cite{LW04}, any smooth solution
of the 3D axisymmetric Navier-Stokes equations must satisfy 
the following compatibility condition at $r=0$:
\begin{equation}
u(0,z,t)=\omega(0,z,t)=\psi(0,z,t)=0.
\label{eqn-bc}
\end{equation}
Moreover, all the even order derivatives of $u$, $\omega$ and $\psi$
with respect to $r$ at $r=0$ must vanish. Therefore, 
we expand the solution $u$, $\omega$ and $\psi$ around $r=0$ as follows:
\begin{eqnarray}
&& u(r,z,t) = r u_1(z,t) + \frac{r^3}{3!} u_3(z,t) + \frac{r^5}{5!} u_5(z,t) + \cdots,
\label{eqn-uexp} \\
&& \omega(r,z,t) = r \omega_1(z,t) + \frac{r^3}{3!} \omega_3(z,t) +
\frac{r^5}{5!} \omega_5(z,t) + \cdots,
\label{eqn-wexp} \\
&& \psi(r,z,t) = r \psi_1(z,t) + \frac{r^3}{3!} \psi_3(z,t) +
\frac{r^5}{5!} \psi_5(z,t) + \cdots .
\label{eqn-pexp}
\end{eqnarray}

Substituting the above expansions into (\ref{eqn-u})-(\ref{eqn-psi}), we obtain
to the leading order the following system of equations:
\begin{eqnarray*}
&& r \left( u_1 \right )_t - r \left (\psi_1 \right )_z u_1 + 2 \psi_1 r
\left ( u_1 \right )_z  = \nu \left ( \frac{4}{3} r u_3 + r
( u_1)_{zz} \right ) + r \left (\psi_1 \right )_z u_1
+ O(r^3)
\\
&& r \left( \omega_1 \right )_t + 2 \psi_1 r \left ( \omega_1 \right )_z  =
\nu \left (  \frac{4}{3} r \omega_3 + r
( \omega_1)_{zz} \right ) + 2 r u_1 \left (u_1 \right )_z
+ O(r^3)
\\
&& - \left (  \frac{4}{3} r \psi_3 + r
( \psi_1)_{zz} + O(r^3 )\right ) = r \omega_1 + O(r^3).
\end{eqnarray*}

By canceling $r$ from both sides and neglecting the higher order terms in $r$, we obtain

\begin{eqnarray*}
&& \left( u_1 \right )_t + 2 \psi_1
\left ( u_1 \right )_z  = \nu \left (\frac{4}{3} u_3 +
(u_1)_{zz} \right ) + 2 \left (\psi_1 \right )_z u_1,
\\
&& \left( \omega_1 \right )_t + 2 \psi_1 \left ( \omega_1 \right )_z  =
\nu \left ( \frac{4}{3} \omega_3 +
(\omega_1)_{zz} \right ) +  \left (u_1^2 \right )_z ,
\\
&& - \left ( \frac{4}{3} \psi_3 +
(\psi_1)_{zz} )\right ) = \omega_1.
\end{eqnarray*}

Note that $u_3 = u_{rrr}(0,z,t)$, $(u_1)_{zz} = u_{rzz}(0,z,t)$. 
If we further make the assumption that the second partial derivative
of $u_1$, $\omega_1$, $\psi_1$ with respect to $z$ is much larger than
the second partial derivative of these functions with respect to $r$,
then we can ignore the coupling in the Laplacian operator to $u_3$,
$\omega_3$ and $\psi_3$ in the above equations. Thus, we obtain our
1D model as follows:
\begin{eqnarray}
&&\left( u_1 \right )_t + 2 \psi_1
\left ( u_1 \right )_z  = \nu (u_1)_{zz} +
2 \left (\psi_1 \right )_z u_1 ,\label{eqn-u10} \\
&&\left( \omega_1 \right )_t + 2 \psi_1 \left ( \omega_1 \right )_z  =
\nu (\omega_1)_{zz} + \left (u_1^2 \right )_z ,\label{eqn-w10}\\
&& - (\psi_1)_{zz} = \omega_1. \label{eqn-psi10}
\end{eqnarray}
We remark that the above assumption implies that the solution 
has an anisotropic scaling, i.e. the solution is more singular 
along the $z$-direction than along the $r$-direction. A possible
scenario is that the solution has a pancake like structure 
perpendicular to the $z$-axis.

Let $\tilde{u} = u_1$, $\tilde{v} = -(\psi_1)_z$, $\tilde{\omega}= \omega_1$,
and $\tilde{\psi} = \psi_1$. By integrating the $\omega_1$ equation with
respect to $z$ and using the relationship
$- \frac{\partial^2}{\partial z^2} \psi_1 = \omega_1$, we can obtain an
evolution equation for $\tilde{v}$. Now the complete set of evolution
equations for $\tilde{u}$, $\tilde{v}$, and $\tilde{\omega}$ are
given by
\begin{eqnarray}
&& \left( \tilde{u} \right )_t + 2 \tilde{\psi} ( \tilde{u})_z  =
\nu  (\tilde{u})_{zz}
- 2 \tilde{v} \tilde{u},
\label{eqn-u1}
\\
&& \left( \tilde{\omega}\right )_t + 2 \tilde{\psi} ( \tilde{\omega} )_z  =
\nu (\tilde{\omega})_{zz} + \left ( \tilde{u}^2\right )_z,
\label{eqn-w1}
\\
&& \left( \tilde{v}\right )_t + 2 \tilde{\psi} ( \tilde{v} )_z  =
\nu (\tilde{v})_{zz} + (\tilde{u})^2 - (\tilde{v})^2
+c(t),
\label{eqn-v1}\\
&& - (\tilde{\psi})_{zz} = \tilde{\omega},
\label{eqn-p1}
\end{eqnarray}
where the constant $c(t)$ is an integration constant which is determined
by enforcing the mean of ${\tilde v} $ equal to zero. For example,
if ${\tilde{\psi}}$ is periodic with period 1 in $z$, then $c(t)$
is given by
\begin{equation}
c(t) = 3 \int_0^1 {\tilde v}^2 dz - \int_0^1 {\tilde u}^2 dz.
\label{eqn-constant}
\end{equation}
Note that the equation for $\tilde{\omega}$ is equivalent to that
for $\tilde{v}$. So it is sufficient to consider the coupled system
for $\tilde{u}$, $\tilde{v}$:
\begin{eqnarray}
&& \left( \tilde{u} \right )_t + 2 \tilde{\psi} ( \tilde{u})_z  =
\nu (\tilde{u})_{zz}
- 2 \tilde{v} \tilde{u}
\label{eqn-u2} \\
&& \left( \tilde{v}\right )_t + 2 \tilde{\psi} ( \tilde{v} )_z  =
\nu (\tilde{v})_{zz} + (\tilde{u})^2 - (\tilde{v})^2
+c(t),
\label{eqn-v2}
\end{eqnarray}
where $\tilde{\psi}$ is related to $\tilde{v}$ by 
$\tilde{v} = -(\tilde {\psi})_z$. By (\ref{eqn-p1}),
we have $\tilde{v}_z = \tilde{\omega}$.

A surprising result is that one can use the above 1D model to construct a
family of {\it exact solutions} for the 3D axisymmetric Navier-Stokes equations.
This is described by the following theorem, which can be verified directly
by substituting (\ref{3D-sol}) into the 3D axisymmetric Navier-Stokes
equations and using the model equation (\ref{eqn-u10})-(\ref{eqn-psi10}).

\vspace{0.2in}
\noindent
{\bf Theorem 1.} {\it Let $u_1$, $\psi_1$ and $\omega_1$ be the solution of the
1D model (\ref{eqn-u10})-(\ref{eqn-psi10}) and define
\begin{equation}
u (r,z,t) = r u_1(z,t), \quad
\omega (r,z,t) = r \omega_1(z,t), \quad
\psi (r,z,t) = r \psi_1(z,t).
\label{3D-sol}
\end{equation}
Then $(u(r,z,t), \;\omega (r,z,t), \; \psi(r,z,t) )$ is
an exact solution of the 3D Navier-Stokes equations.}

\vspace{0.2in}
Theorem 1 tells us that the 1D model (\ref{eqn-u10})-(\ref{eqn-psi10})
preserves some essential nonlinear structure of the original 3D 
axisymmetric Navier-Stokes equations. As we will see later, the 
nonlinear structure of the 1D model plays a critical role in 
stabilizing the solution for large times, although the same nonlinearity
can lead to large dynamic growth for early times.

\section{Properties of the Model Equation}

In this section, we will study some properties of the 1D model equations.
We will first consider the properties of some further simplified models obtained from these equations.
Both numerical and analytical studies are presented for these simplified
models. Based on the understanding of the simplified models, we prove the
global existence of the inviscid Lagrangian model, which sheds useful
light into our global existence analysis for the full 1D model with or
without viscosity.

\subsection{The ODE model}

To start with, we consider an ODE model by ignoring the convection
and diffusion terms.
\begin{eqnarray}
&& \left( \tilde{u} \right )_t = - 2 \tilde{v} \tilde{u} \label{u-ode} \\
&& \left( \tilde{v}\right )_t = (\tilde{u})^2 - (\tilde{v})^2, \label{v-ode}
\end{eqnarray}
with initial condition ${\tilde u}(0)={\tilde u}_0$ and
${\tilde v}(0)={\tilde v}_0$.

Clearly, if ${\tilde u}_0 = 0$, then ${\tilde u}(t) = 0 $ for all $t>0$.
In this case, the equation for ${\tilde v}$ is decoupled from ${\tilde u}$
completely, and will blow up in finite time if ${\tilde v}_0 < 0$.
In fact, if ${\tilde v}_0 < 0$ and ${\tilde u}_0$ is very small, then
the solution can experience very large growth dynamically. The growth
can be made arbitrarily large if we choose ${\tilde u}_0$ to be arbitrarily
small. However, the special nonlinear structure of the ODE system
has an interesting cancellation property which has a stabilizing 
effect of the solution for large times. This is described by the
following theorem.

\vspace{0.2in}
\noindent
{\bf Theorem 2.} Assume that ${\tilde u}_0 \neq 0$. Then the solution
$({\tilde u}(t), {\tilde v}(t))$ of the ODE system (\ref{u-ode})-(\ref{v-ode})
exists for all times. Moreover, we have
\begin{equation}
\lim_{t \rightarrow +\infty} {\tilde u} (t) = 0, \quad
\lim_{t \rightarrow +\infty} {\tilde v} (t) = 0 .
\end{equation}

\vspace{0.1in}
\noindent
{\bf Proof.} There are several ways to prove this theorem. The simplest
way is to reformulate the problem in terms of complex
variables\footnote{We thank Prof. Tai-Ping Liu
for suggesting the use of complex variables}. Let
\[
w = {\tilde u} + i {\tilde v}.
\]
Then the ODE system (\ref{u-ode})-(\ref{v-ode}) is reduced to
the following complex nonlinear ODE:
\begin{equation}
\frac{ d w}{dt } = i w^2 , \quad w(0)=w_0,
\label{ode_complex}
\end{equation}
which can be solved analytically. 
The solution has the form
\begin{equation}
w(t) = \frac{w_0}{1-iw_0 t}.
\end{equation}
In terms of the original variables, we have
\begin{eqnarray}
{\tilde u}(t) = \frac{{\tilde u}_0(1+{\tilde v}_0t) - {\tilde u}_0{\tilde v}_0t}
{(1+{\tilde v}_0t)^2 + ({\tilde u}_0t)^2} ,\label{u-exact} \\
{\tilde v}(t) = \frac{{\tilde v}_0(1+{\tilde v}_0t) + {\tilde u}_0^2 t}
{(1+{\tilde v}_0t)^2 + ({\tilde u}_0t)^2} .\label{v-exact}
\end{eqnarray}

It is clear from (\ref{u-exact})-(\ref{v-exact}) that the solution of the
ODE system (\ref{u-ode})-(\ref{v-ode}) exists for all times and
decays to zero as $t \rightarrow +\infty$ as long as ${\tilde u}_0 \neq 0$.
This completes the proof of Theorem 2.

\vspace{0.2in}
\noindent
{\bf Remark 1}. 
Note that the ODE model (\ref{ode_complex}) has some similarity 
with the Constantin-Lax-Majda model \cite{CLM85}, which has the form
$u_t = u H(u)$, where $H$ is the Hilbert transform. By letting
$w= H(u) + iu$ and using the property of the Hilbert transform, 
Constantin-Lax-Majda show that their model can be written as
the imaginary part of the complex ODE: $w_t = \frac{1}{2} w^2$.
It is interesting to note that both models ignore the convection
term and they have solutions that blow up at a finite time for
initial condition satisfying $u(z_0) = 0$ and $H(u)(z_0)>0$
for some $z_0$. However, as we will show later, the convection
term plays an important role in stabilizing the 1D model
and should not be neglected in our study of the Euler equations. 
By including the convection term in the 1D model, we will show 
in section 3.3 ans section 4 that no finite time blow-up 
can occur from smooth initial data.

As we can see from (\ref{u-exact})-(\ref{v-exact}),
the solution can grow very fast in a very short time if 
$\tilde u_0$ is small, but $\tilde v_0$ is large and negative.
For example, if we
let ${\tilde v}_0 = - 1/\epsilon$ and ${\tilde u}_0 =\epsilon$ for
$\epsilon > 0$ small, we obtain at $t=\epsilon$ 
\[
{\tilde u}(\epsilon ) = 1/\epsilon^3, \quad
{\tilde v}(\epsilon ) = 1/\epsilon .
\]
We can see that within $\epsilon $ time, ${\tilde u}$ grows 
from its initial value of order $\epsilon$ to $O(\epsilon^{-3})$, a factor of $\epsilon^{-4}$ amplification. 

\vspace{0.2in}
\noindent
{\bf Remark 2}. The key ingredient in obtaining the global existence in
Theorem 2 is that the coefficient on the right hand side of (\ref{u-ode})
is less than $-1$. For this ODE system, there are two distinguished
phases. In the first phase, if ${\tilde v}$ is negative and large
in magnitude, but ${\tilde u}$ is small, then ${\tilde v}$ can
experience tremendous dynamic growth, which is essentially governed by
\[
{\tilde v}_t = - {\tilde v}^2.
\]
However, as ${\tilde v}$ becomes very large and negative, it will 
induce a rapid growth in ${\tilde u}$. The nonlinear structure of 
the ODE system is such that ${\tilde u}$ will eventually grow even 
faster than ${\tilde v}$ and force $(\tilde{u})^2 - (\tilde{v})^2 < 0$ 
in the second phase. From this time on, $\tilde{v}$ will increase in 
time and eventually become positive. Once $\tilde{v}$ becomes positive, 
the nonlinear term, $- \tilde{v}^2$, becomes stabilizing for $\tilde{v}$.
Similarly, the nonlinear term, $-2 \tilde u \tilde v$, becomes stabilizing
for $\tilde u$. This subtle dynamic stability property of the ODE
system can be best illustrated by the
phase diagram in Figure \ref{fig.phasediagram}.

In Appendix A, we  prove the same result for
a more general ODE system of the following form:
\begin{eqnarray}
&& \left( \tilde{u} \right )_t = - d \tilde{v} \tilde{u} \label{u-ode1} \\
&& \left( \tilde{v}\right )_t = (\tilde{u})^2 - (\tilde{v})^2, \label{v-ode1}
\end{eqnarray}
for any constant $d \geq 1$. However, if $d < 1$, it is possible to
construct a family of solutions for the ODE
systems (\ref{u-ode1})-(\ref{v-ode1}) which blow up in a finite time.

\begin{figure}
\begin{center}
\includegraphics[width=10cm]{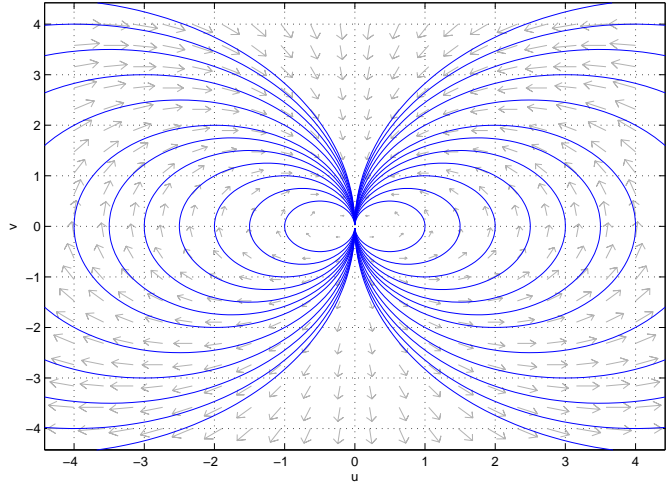}
\end{center}
\caption{The Phase Diagram for the ODE system.
\label{fig.phasediagram}}
\end{figure}

\subsection{The Reaction Diffusion Model}

In this subsection, we consider the reaction-diffusion system:
\begin{eqnarray}
&& \left( \tilde{u} \right )_t = \nu \tilde{u}_{zz}
- 2 \tilde{v} \tilde{u}, \label{eqn-unew} \\
&& \left( \tilde{v}\right )_t = \nu \tilde{v}_{zz} +
(\tilde{u})^2 - (\tilde{v})^2. \label{eqn-vnew}
\end{eqnarray}

As we can see for the corresponding ODE system, the structure of the
nonlinearity plays an essential role in obtaining global existence.
Intuitively, one may think that the diffusion term would
help to stabilize the dynamic growth induced by the nonlinear terms. 
However, because the nonlinear ODE system in
the absence of viscosity is very unstable, the diffusion term 
can actually have a destabilizing effect. Below we will 
demonstrate this somewhat surprising fact through careful
numerical experiments.

\begin{figure}
\begin{center}
\includegraphics[width=8cm]{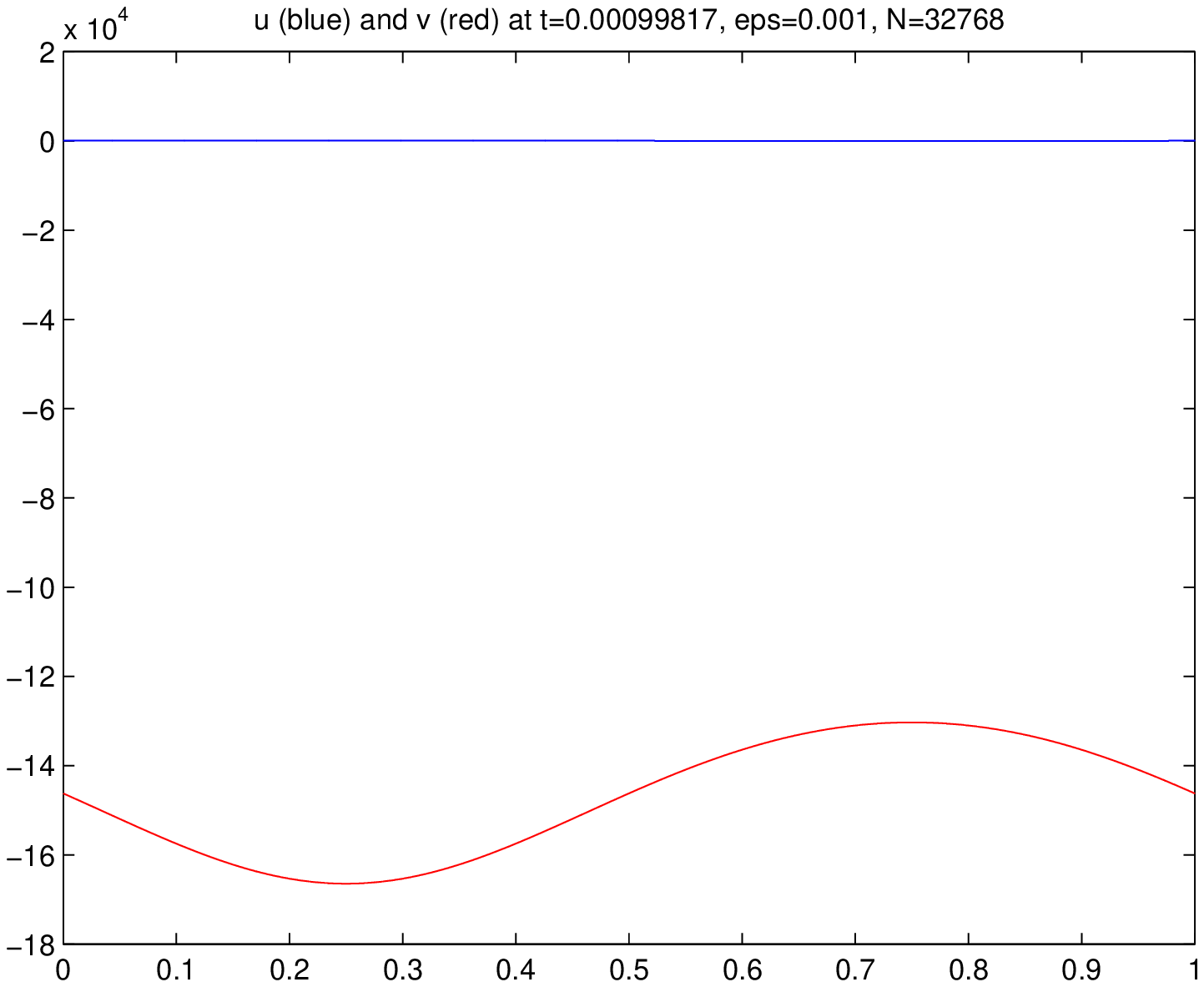}
\includegraphics[width=8cm]{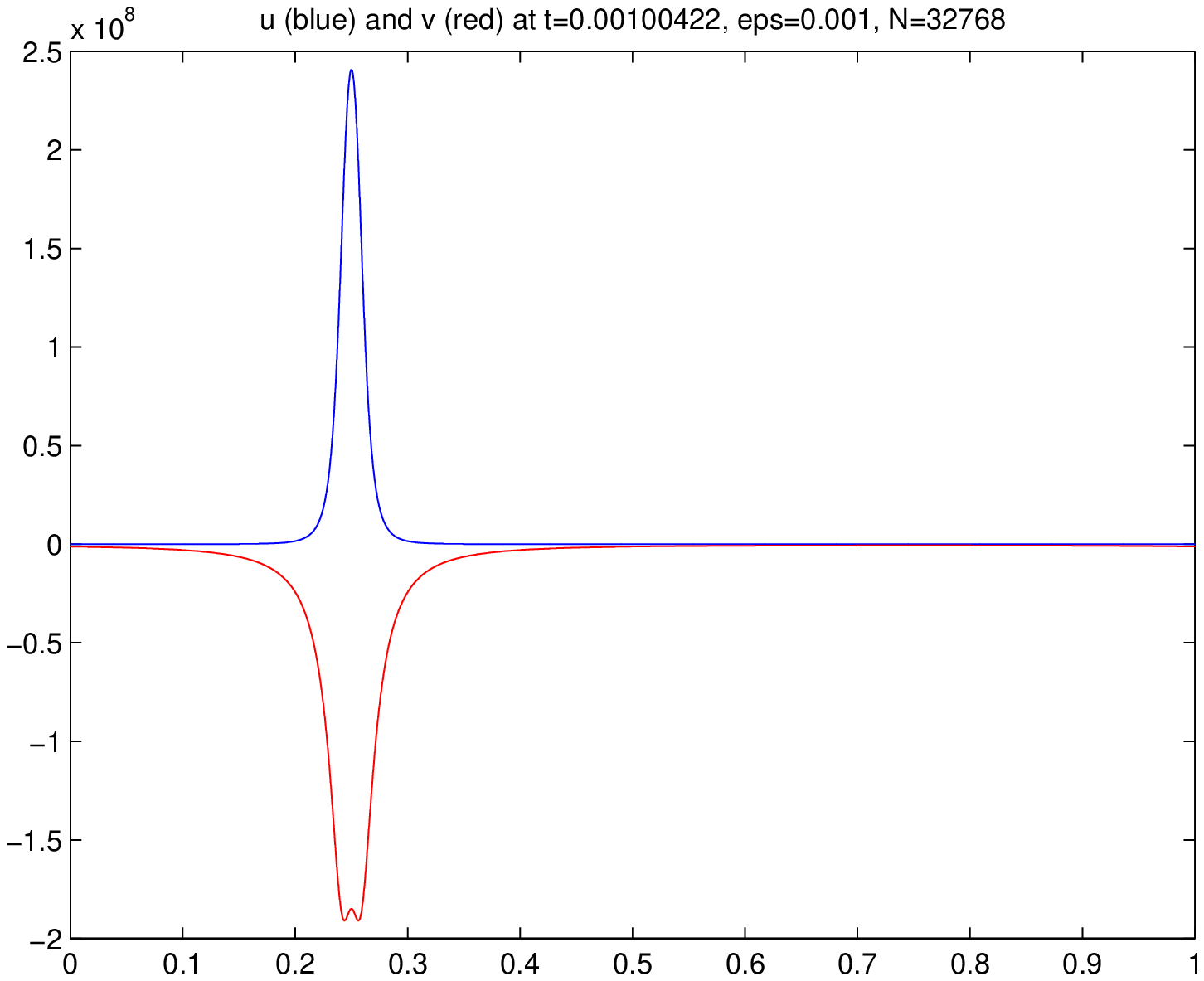}
\end{center}
\caption{The solutions at $t=0.00099817$, and $t=0.0010042$, $N=32768$, $\nu = 1$.
\label{fig.uv1}}
\end{figure}

\begin{figure}
\begin{center}
\includegraphics[width=8cm]{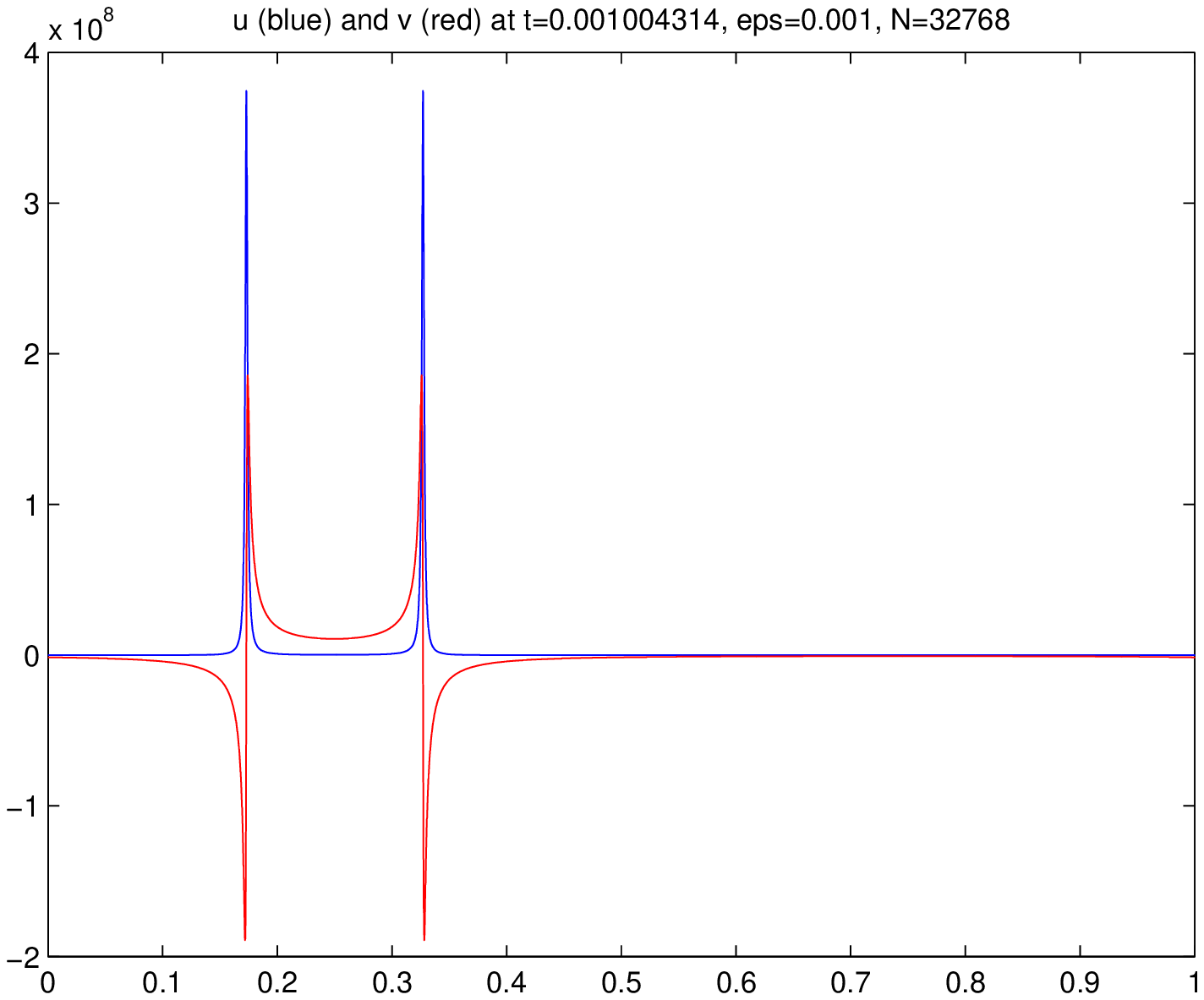}
\includegraphics[width=8cm]{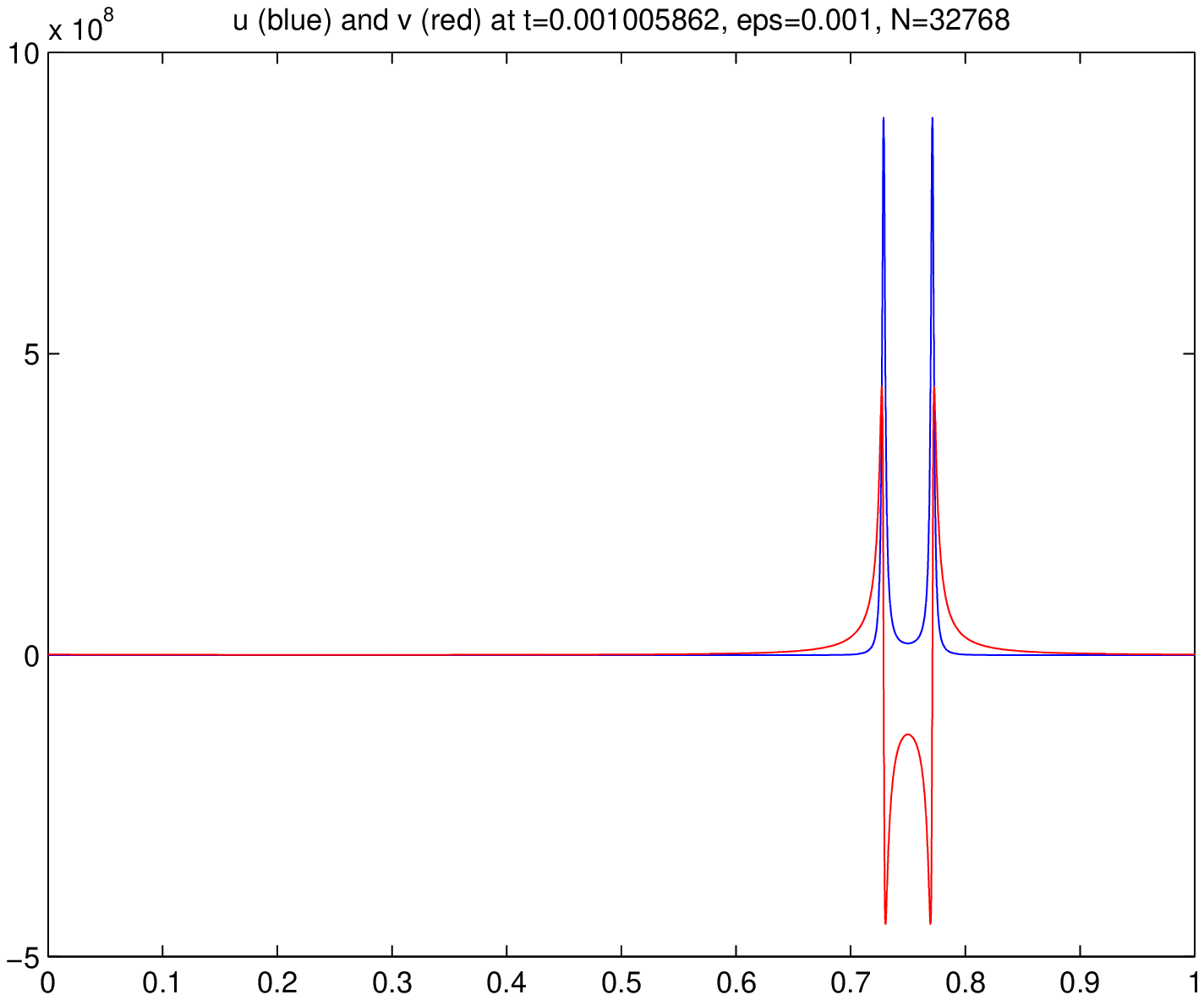}
\end{center}
\caption{The solutions at $t=0.001004314$ and $t=0.001005862$, $N=32768$, $\nu = 1$.
\label{fig.uv2}}
\end{figure}

\begin{figure}
\begin{center}
\includegraphics[width=8cm]{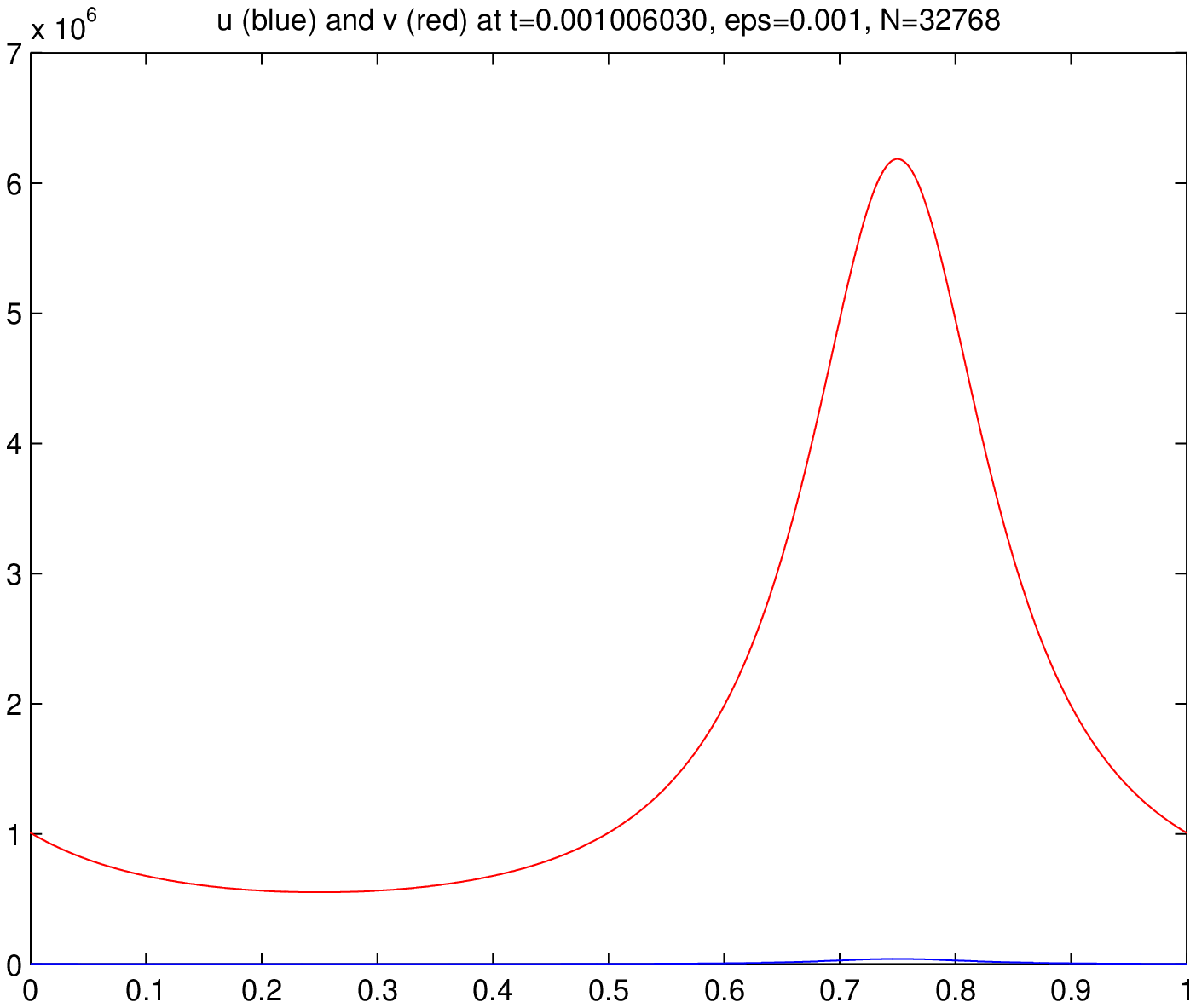}
\includegraphics[width=8cm]{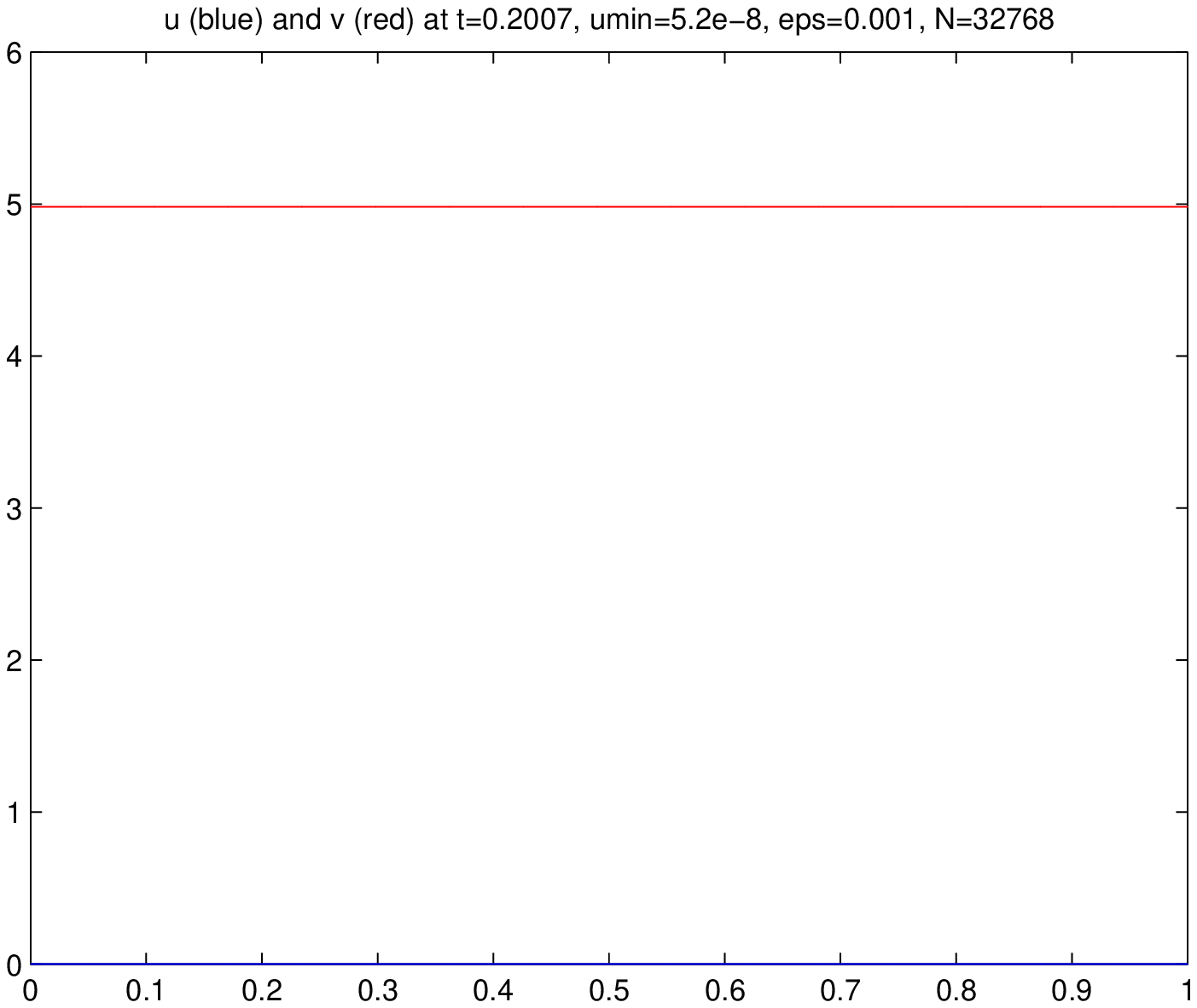}
\end{center}
\caption{The solutions at $t=0.00100603$ and $t=0.2007$, $N=32768$, $\nu = 1$.
\label{fig.uv3}}
\end{figure}

In Figures (\ref{fig.uv1})-(\ref{fig.uv3}), we plot a time
sequence of solutions for the above reaction diffusion system
with the following initial data
\[
{\tilde u}_0 (z) = \epsilon (2+\sin (2 \pi z)), \quad
{\tilde v}_0 (z) = -\frac{1}{\epsilon} - \sin (2 \pi z),
\]
where $\epsilon = 0.001$.
For this initial condition, the solution is periodic in $z$ with
period one. We use a pseudo-spectral
method to discretize the coupled system (\ref{eqn-unew})-(\ref{eqn-vnew})
in space and use the simple forward Euler discretization for the
nonlinear terms and the backward Euler discretization for the diffusion
term. In order to resolve the nearly singular solution structure, we
use $N=32,768$ grid points with an adaptive time step satisfying
\[
\Delta t_n \left ( | \max\{ {\tilde u}^n\} | + |\min \{ {\tilde u}^n\}|
+ | \max\{ {\tilde v}^n\} | + |\min \{ {\tilde v}^n\}|
\right ) \leq 0.01,
\]
where ${\tilde u}^n$ and ${\tilde v}^n$ are the numerical solution
at time $t_n$ and $t_n = t_{n-1} + \Delta t_{n-1}$ with the
initial time stepsize $\Delta t_0 = 0.01 \epsilon$. During
the time iterations, the smallest time step is as small
as $O(10^{-10})$. 

From Figure \ref{fig.uv1}, we can see that the magnitude of the
solution $\tilde v$ increases rapidly by a factor of 150 within
a very short time ($t=0.00099817$). As the solution $\tilde v$
becomes large and negative, the solution $\tilde u$ increases
much more rapidly than $\tilde{v}$. By time $t=0.0010042$, $\tilde u$ has increased
to about $2.5\times 10^8$ from its initial condition which is of
magnitude $10^{-3}$. This is a factor of $2.5 \times 10^{11}$ increase.
At this time, the minimum of $\tilde v$ has reached 
$-2\times 10^8$. Note that since $\tilde u$ has outgrown $\tilde v$
in magnitude, the nonlinear term, ${\tilde u}^2 - {\tilde v}^2$,
on the right hand side of
the $\tilde v$-equation has changed sign. This causes the
solution $\tilde v$ to split.
By the time $t=0.001004314$ (see Figure \ref{fig.uv2}), 
both $\tilde u$ and $\tilde v$ have 
split and settled down to two relatively stable traveling wave 
solutions. The wave on the left will travel to the left while 
the wave on the right will travel to the right. Due to the 
periodicity in $z$, the two traveling waves approach each 
other from the right side of the domain. The ``collision''
of these two traveling waves tends to annihilate each other. In
particular, the negative part of $\tilde v$ is effectively
eliminated during this nonlinear interaction. By the time
$t=0.00100603$ (see Figure \ref{fig.uv3}), 
the solution $\tilde v$ becomes all positive.
Once $\tilde v$ becomes positive, the effect of nonlinearity
becomes stabilizing for both $\tilde u$ and $\tilde v$,
as in the case of the ODE system. From then on, the solution 
decays rapidly. By $t=0.2007$, the magnitude of $\tilde u$ 
is as small as $5.2 \times 10^{-8}$, and $\tilde v$ becomes almost 
a constant function with value close to 5. From this time on, 
$\tilde u$ is essentially decoupled from $\tilde v$ and will 
decay like $O(1/t)$.

\subsection{The Lagrangian Convection Model}

Next, we consider the 1D model equations in the absence of viscosity.
The corresponding equations are given as follows:
\begin{eqnarray}
&& \tilde{u}_t + 2 \tilde \psi {\tilde u}_z =
- 2 \tilde{v} \tilde{u} \label{eqn-ueuler} \\
&& \tilde{v}_t + 2 \tilde \psi {\tilde v}_z =
\tilde{u}^2 - \tilde{v}^2 + c(t), \label{eqn-veuler}
\end{eqnarray}
where $\tilde{v} = - \tilde \psi_z$, and 
$c(t)$ is defined in (\ref{eqn-constant}) 
to ensure that $\int_0^1 \tilde{v} dz = 0$.

Introduce the Lagrangian flow map
\begin{eqnarray}
&& \frac {\partial z(\alpha ,t) }{\partial t} = 2 \tilde \psi (z(\alpha ,t),t),
\label{eqn-flowmap}\\
&& z(\alpha ,0)=\alpha.
\end{eqnarray}
Differentiating (\ref{eqn-flowmap}) with respect to $\alpha$, we get
\[
\frac{\partial z_\alpha }{\partial t} = 2 z_\alpha
\frac{\partial \tilde \psi }{\partial z} (z(\alpha ,t),t)
= -2 z_\alpha \tilde v (z(\alpha ,t),t) .
\]

Denote $v (\alpha ,t)= \tilde v (z(\alpha ,t),t)$,
$u (\alpha ,t)= \tilde u (z(\alpha ,t),t)$, and
$J(\alpha ,t)= z_\alpha (\alpha ,t)$. Then we can show that
$J$, $u$, and $v$ satisfy the following system of equations:
\begin{eqnarray}
&& \frac {\partial J(\alpha ,t) }{\partial t} = -2 J(\alpha ,t) v(\alpha ,t),
\label{eqn-J}\\
&& \frac {\partial u(\alpha ,t) }{\partial t} = -2 u(\alpha ,t) v(\alpha ,t),
\label{eqn-uL}\\
&& \frac {\partial v(\alpha ,t) }{\partial t} = u^2 -  v^2 + 3 \int_0^1 v^2 J d \alpha
- \int_0^1 u^2 J d \alpha ,
\label{eqn-vL}
\end{eqnarray}
with initial data $J(\alpha ,0) = 1$, $u(\alpha ,0) = {\tilde u}_0 (\alpha )$
and $v(\alpha ,0) = {\tilde v}_0 (\alpha )$. 
Since $\int_0^1 \tilde v (z, t) d z = 0$, we have
\begin{equation}
\int_0^1 v(\alpha ,t) J(\alpha ,t) d \alpha = 0,
\label{eqn-newmean}
\end{equation}
which implies that
\begin{equation}
\int_0^1 J(\alpha ,t) d \alpha \equiv  \int_0^1 J(\alpha ,0) d \alpha = 1.
\label{eqn-Jmean}
\end{equation}
It is interesting to note that the 1D model formulated in the Lagrangian 
coordinate retains some of the essential properties of the ODE system.
In the following, we will explore the special nonlinear structure
of the model equation to prove the global well-posedness of the
1D model in the Lagrangian form. As we will see, the understanding of
the 1D model in the Lagrangian form gives critical insight in our 
understanding of the full 1D model.

\vspace{0.2in}
\noindent
{\bf Theorem 3.} {\it Assume that
$\tilde u (z,0)$ and $\tilde v (z,0)$ are in 
$C^m[0,1]$ with $m \geq 1$ and periodic with period 1.
Then the solution $(\tilde u, \tilde v)$ of the 1D inviscid model 
will be in $C^m [0,1]$ for all times. }

\vspace{0.2in}
\noindent
{\bf Proof.} Differentiating the $\tilde{u}$ and $\tilde v$ equations
with respect to $\alpha$, we get
\begin{eqnarray}
&&\frac{d {\tilde u}_\alpha}{dt} = - 2 \tilde v {\tilde u}_\alpha -
2 {\tilde u} {\tilde v}_\alpha ,\label{eqn-ua}\\
&&\frac{d {\tilde v}_\alpha}{dt} =  2 \tilde u {\tilde u}_\alpha -
2 {\tilde v} {\tilde v}_\alpha \label{eqn-va} .
\end{eqnarray}
Multiplying (\ref{eqn-ua}) by ${\tilde u}_\alpha$ and
(\ref{eqn-va}) by ${\tilde v}_\alpha$, and adding the resulting
equations, we have
\begin{eqnarray}
\frac{1}{2} \frac{d}{dt} \left ( {\tilde u}_\alpha^2
+ {\tilde v}_\alpha^2 \right ) = -2 {\tilde v}
\left ( {\tilde u}_\alpha^2 + {\tilde v}_\alpha^2 \right ).
\end{eqnarray}
Therefore, we obtain
\begin{equation}
\frac{1}{2} \frac{d}{dt} \log \left ( {\tilde u}_\alpha^2
+ {\tilde v}_\alpha^2 \right ) = - 2 {\tilde v} .
\end{equation}
Integrating from 0 to $t$, we get
\begin{eqnarray}
\left (\sqrt{{\tilde u}_\alpha^2 + {\tilde v}_\alpha^2 } \right ) (\alpha, t)
= \left (\sqrt{({\tilde u}_0)_\alpha^2 + ({\tilde v}_0)_\alpha^2 }
\right ) e^{-2 \int_0^t {\tilde v} (\alpha , s ) ds}
= \left (\sqrt{({\tilde u}_0)_\alpha^2 + ({\tilde v}_0)_\alpha^2 }
\right ) J(\alpha, t) ,
\label{uv_alpha}
\end{eqnarray}
where we have used
\[
J(\alpha, t) = e^{-2 \int_0^t {\tilde v} (\alpha , s ) ds},
\]
which follows from (\ref{eqn-J}) and $J(\alpha,0) \equiv 1$.
Using (\ref{eqn-Jmean}), we further obtain
\begin{equation}
\int_0^1 \sqrt{{\tilde u}_\alpha^2 + {\tilde v}_\alpha^2 } \; d \alpha
\leq \|\sqrt{({\tilde u}_0)_\alpha^2 + ({\tilde v}_0)_\alpha^2 }\;\|_{L^\infty}
\int_0^1 J(\alpha ,t) d \alpha =
\|\sqrt{({\tilde u}_0)_\alpha^2 + ({\tilde v}_0)_\alpha^2 } \;\|_{L^\infty}.
\label{bound-uv}
\end{equation}
In particular, we have
\begin{equation}
\int_0^1 | {\tilde v}_\alpha | d \alpha \leq \int_0^1
\sqrt{{\tilde u}_\alpha^2 + {\tilde v}_\alpha^2 } d \alpha
\leq \|\sqrt{({\tilde u}_0)_\alpha^2 + ({\tilde v}_0)_\alpha^2 } \;\|_{L^\infty}.
\label{bound-v-tv}
\end{equation}
Since $\int_0^1 {\tilde v} J d \alpha = 0$ and $J >0$, there exists
$\alpha_0(t) \in [0,1]$ such that ${\tilde v} (\alpha_0(t),t) = 0$. Therefore,
we get
\begin{eqnarray}
| {\tilde v} (\alpha ,t) |  =  | \int_{\alpha_0}^\alpha {\tilde v}_\alpha
d \alpha' |
\leq \int_0^1 | {\tilde v} (\alpha',t) | d \alpha' \leq
\|\sqrt{({\tilde u}_0)_\alpha^2 + ({\tilde v}_0)_\alpha^2 }\;\|_{L^\infty}.
\end{eqnarray}
This proves that
\begin{equation}
\| {\tilde v} \|_{L^\infty} \leq
\|\sqrt{({\tilde u}_0)_\alpha^2 + ({\tilde v}_0)_\alpha^2 }\;\|_{L^\infty}.
\end{equation}
Using the equations for $J$ and $\tilde{u}$, we also obtain
\begin{eqnarray}
&& e^{-2 t C_0} \leq J (\alpha ,t) \leq e^{2 t C_0},  \\
&&\| \tilde u \|_{L^\infty} \leq \| {\tilde u}_0 \; \|_{L^\infty}
e^{2 t C_0},
\end{eqnarray}
where
$C_0=\|\sqrt{({\tilde u}_0)_\alpha^2 + ({\tilde v}_0)_\alpha^2 }\;\|_{L^\infty}.$

The bound on $J(\alpha, t)$ in turn gives bound on
${\tilde u}_\alpha^2 + {\tilde v}_\alpha^2 $ through (\ref{uv_alpha}).
We can then bootstrap to obtain regularity of the solution in 
higher order norms.  This completes the proof of Theorem 3.

\begin{figure}
\begin{center}
\includegraphics[width=8cm]{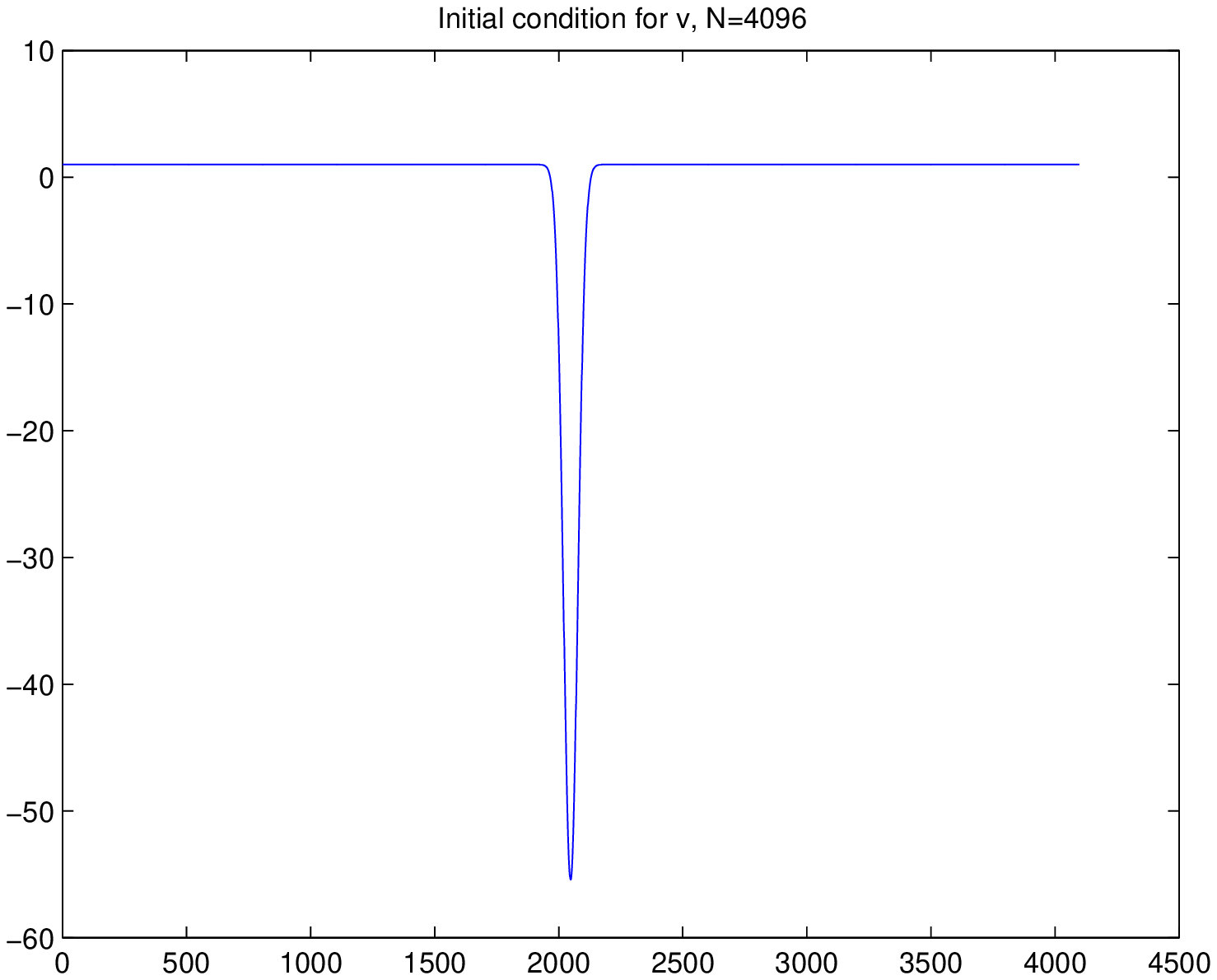}
\includegraphics[width=8cm]{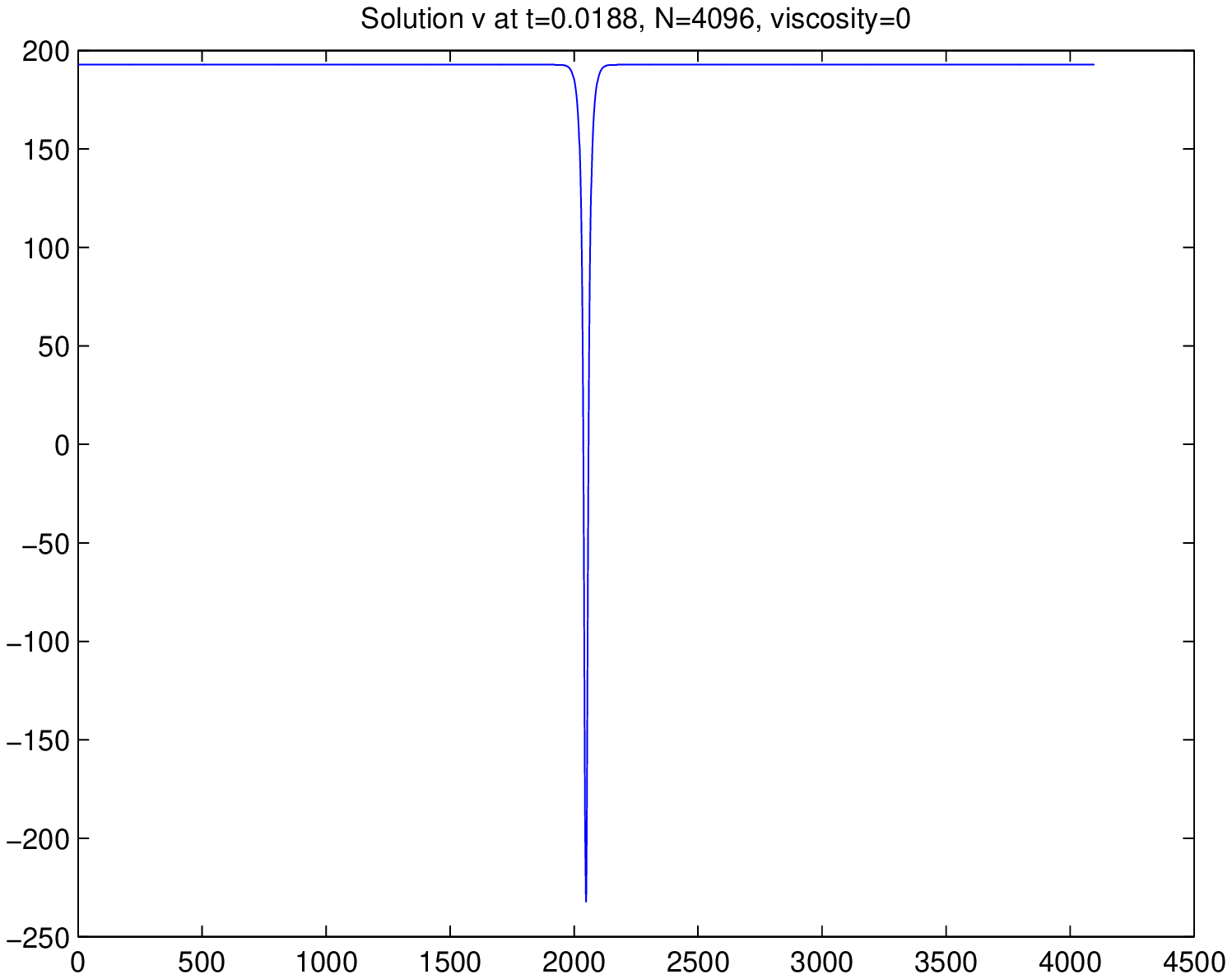}
\end{center}
\caption{The Lagrangian solution at $t=0$ and $t=0.0188$, $N=4096$.
\label{fig.uv5}}
\end{figure}

\begin{figure}
\begin{center}
\includegraphics[width=8cm]{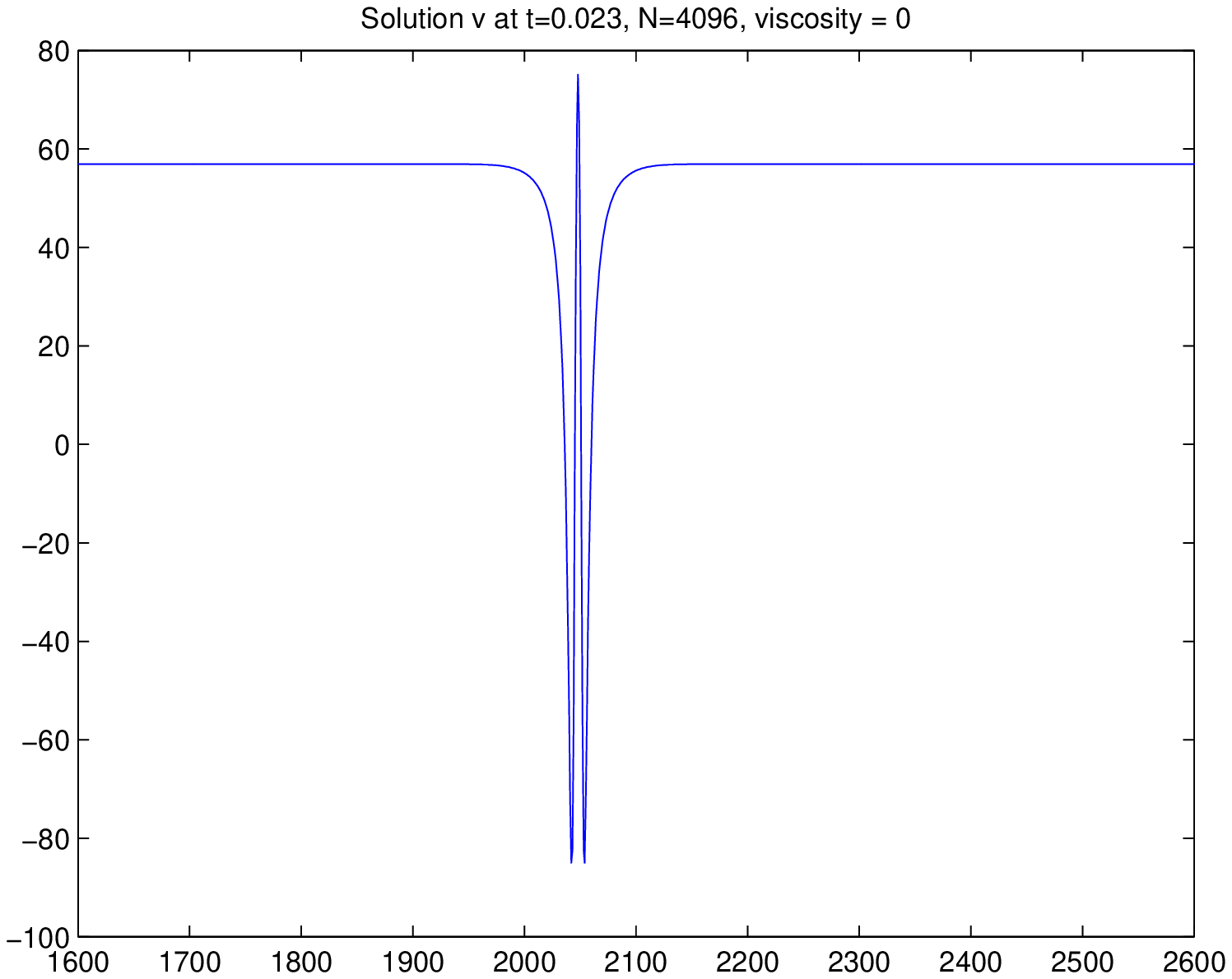}
\includegraphics[width=8cm]{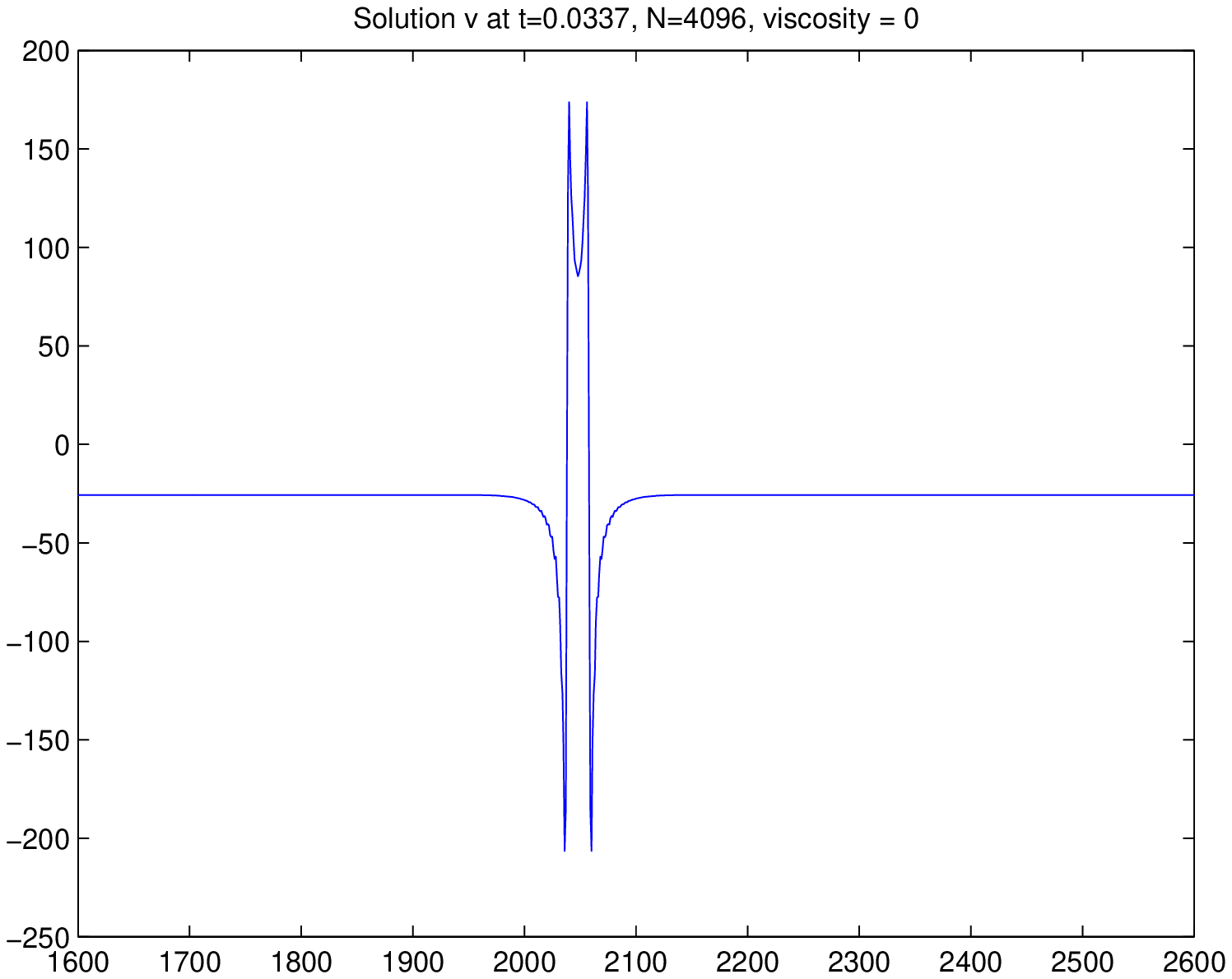}
\end{center}
\caption{The Lagrangian solution at $t=0.023$ and $t=0.0337$, $N=4096$.
\label{fig.uv6}}
\end{figure}

\begin{figure}
\begin{center}
\includegraphics[width=8cm]{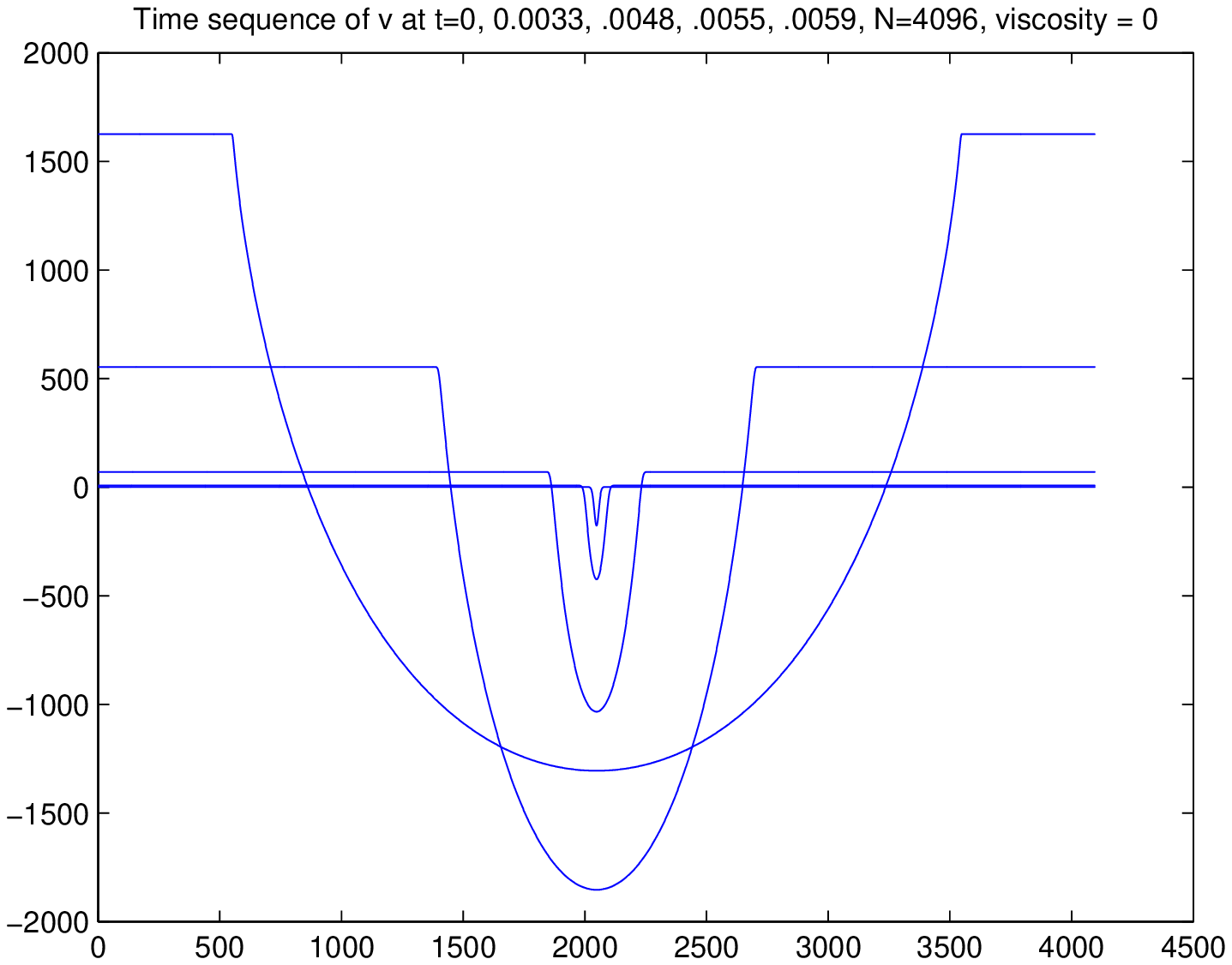}
\includegraphics[width=8cm]{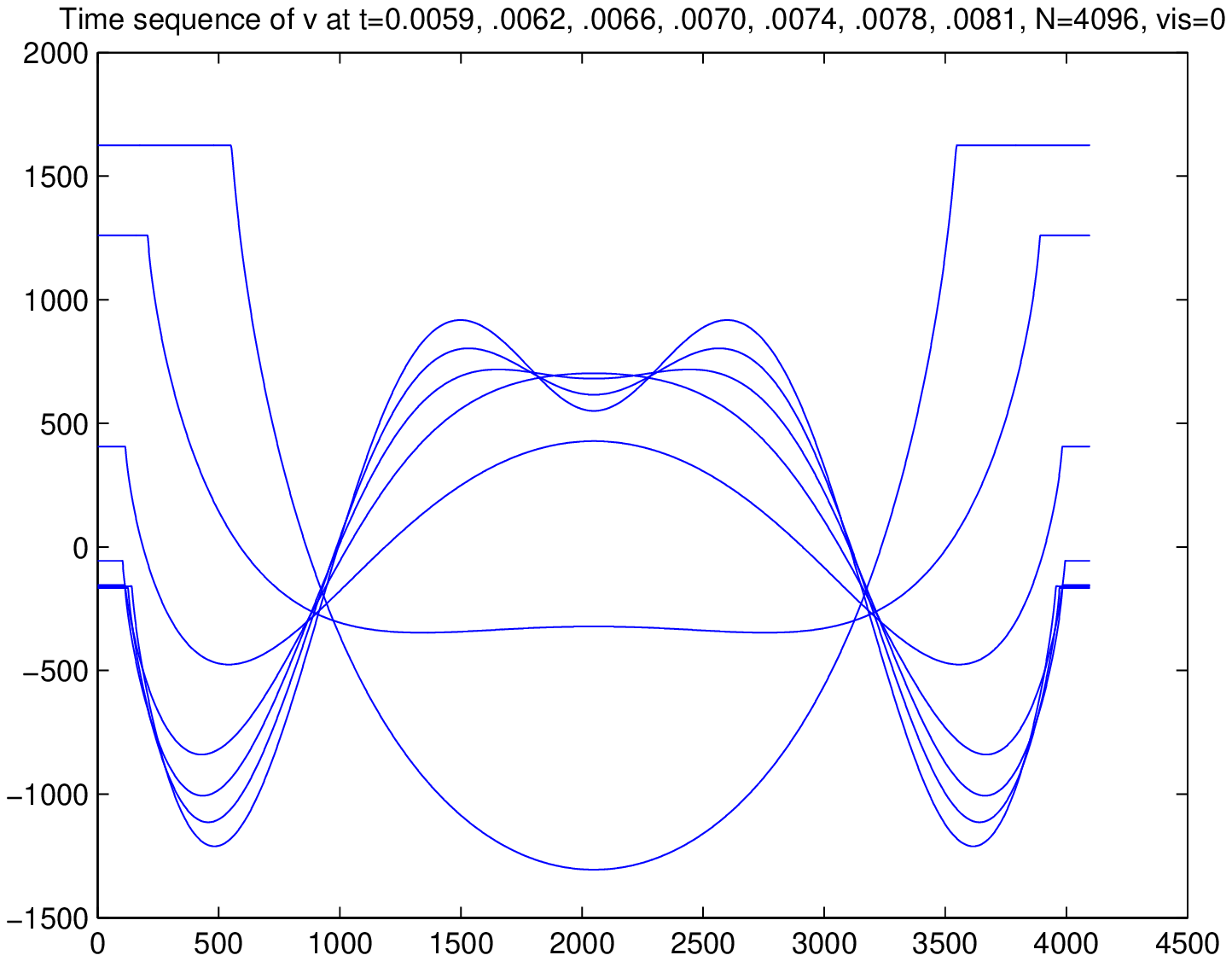}
\end{center}
\caption{The Eulerian solution,  $N=4096$.
\label{fig.uv7}}
\end{figure}

\begin{figure}
\begin{center}
\includegraphics[width=12cm]{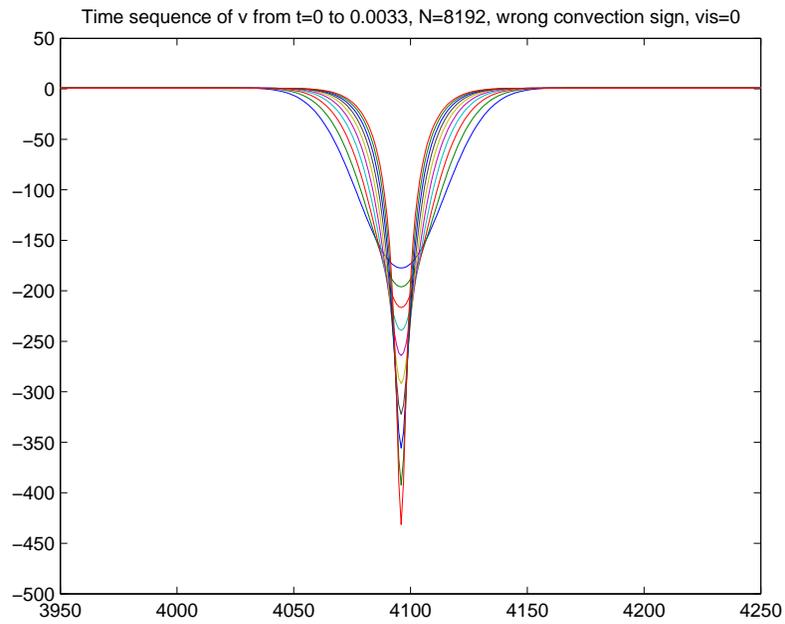}
\end{center}
\caption{The Lagrangian solutions with the wrong sign, $N=4096$.
\label{fig.uv8}}
\end{figure}

Next, we illustrate the behavior of the solution through  numerical
computations. We use a pseudo-spectral method to discretize
in space and a second order Runge-Kutta discretization in time
with an adaptive time-stepping.
In Figures \ref{fig.uv5} and \ref{fig.uv6}, we plot a sequence of snapshots of the
solution for the inviscid model (\ref{eqn-J})-(\ref{eqn-vL}) in the Lagrangian
coordinate using the following initial data
\[
u(\alpha,0)=1, \quad v(\alpha ,0) = 1 - \frac{1}{\delta}\exp^{-(x-0.5)^2/\epsilon} ,
\]
with $\epsilon = 0.0001$ and $\delta = \sqrt{\epsilon \pi}$. 
We can see that the
solution experiences a similar splitting process as in the reaction
diffusion model. In Figure \ref{fig.uv7}, we perform a similar
computation in the Eulerian coordinate with $\epsilon =0.00001$.
We can see that as the solution $\tilde v$ grows large and negative,
the initial sharp profile of $\tilde v$ becomes wider and smoother.
This is a consequence of the incompressibility of the fluid flow.
If we change the sign of the convection velocity from $2\tilde \psi $
to $- 2\tilde \psi $, the profile of $\tilde v$ becomes focused
dynamically and develops an unphysical ``shock-like'' solution, which
seems to evolve into  a finite time blowup, see Figure \ref{fig.uv8}.

%
%

\section{Global Well-Posedness of the 1D Viscous Model}

Based on the understanding we have gained from the previous sections,
we are ready to present a complete proof of the global well-posedness
of the full 1D model. It is not easy to obtain global regularity of the
1D model by using an energy type of estimates. If we multiply
the $\tilde u$-equation by $\tilde u$, and the $\tilde v$-equation
by $\tilde v$, and integrate over $z$, we would arrive at
\begin{eqnarray}
&& \frac{1}{2} \frac{d}{d t}\int_0^1 {\tilde u}^2 dz = - 3 \int_0^1
(\tilde u )^2 \tilde v dz - \nu \int_0^1 {\tilde u}_z^2 dz,\\
&& \frac{1}{2} \frac{d}{d t} \int_0^1 {\tilde v}^2 dz =
\int_0^1 {\tilde u}^2 \tilde v dz - 3 \int_0^1
(\tilde v )^3 dz - \nu \int_0^1 {\tilde v}_z^2 dz.
\end{eqnarray}

Even for this 1D model, the energy estimate shares the some essential
difficulty as the 3D Navier-Stokes equations. It is not clear how
to control the nonlinear vortex stretching like terms by the diffusion
terms. On the other hand, if we assume that
\[
\int_0^T \| \tilde v \|_{L^\infty} dt < \infty ,
\]
similar to the Beale-Kato-Majda non-blowup condition for 
vorticity \cite{BKM84}, then
one can easily show that there is no blow-up up to $t= T$.

In order to obtain the global regularity of the 1D model, we need to
use a local estimate. We will prove that if the initial conditions
for $\tilde u$ and $\tilde v$ are in $C^m $ with $m \ge 1$, then
the solution will remain in $C^m$ for all times.

\vspace{0.2in}
\noindent
{\bf Theorem 4.} {\it Assume that $\tilde u (z,0)$ and $\tilde v (z,0)$ are
in $C^m[0,1]$ with $m \geq 1$ and periodic with period 1. Then the 
solution $(\tilde u, \tilde v)$ of the 1D model will be in $C^m[0,1]$ for all times.}

\vspace{0.2in}
\noindent
{\bf Proof.} Motivated by our analysis for the inviscid Lagrangian
model, we will try to obtain 
{\it a priori} estimate for the nonlinear term
${\tilde u}_z^2 + {\tilde v}_z^2 $.
Differentiating the $\tilde u$-equation
and the $\tilde v$-equation with respect to $z$, we get
\begin{eqnarray}
&& ({\tilde u}_z)_t + 2 \tilde \psi (\tilde u_z )_z - 2 \tilde v
\tilde u_z = - 2 \tilde u_z \tilde v - 2 \tilde u \tilde v_z
+\nu (\tilde u_z)_{zz},\\
&& ({\tilde v}_z)_t + 2 \tilde \psi (\tilde v_z )_z - 2 \tilde v
\tilde v_z =  2 \tilde u \tilde u_z - 2 \tilde v \tilde v_z
+\nu (\tilde v_z)_{zz} .
\end{eqnarray}
Note that one of the nonlinear terms resulting from differentiating
the convection term cancels one of the nonlinear terms on the 
right hand side. After canceling the same nonlinear term from
both sides, we obtain 
\begin{eqnarray}
&& ({\tilde u}_z)_t + 2 \tilde \psi (\tilde u_z )_z
 = - 2 \tilde u \tilde v_z
+\nu (\tilde u_z)_{zz},
\label{eqn-uz}\\
&& ({\tilde v}_z)_t + 2 \tilde \psi (\tilde v_z )_z
 =  2 \tilde u \tilde u_z
+\nu (\tilde v_z)_{zz}.
\label{eqn-vz}
\end{eqnarray}
Multiplying (\ref{eqn-uz}) by $\tilde u_z$ and
(\ref{eqn-vz}) by $\tilde v_z$, we have
\begin{eqnarray}
&& \frac{1}{2}({\tilde u}_z^2)_t + \tilde \psi (\tilde u_z^2 )_z
 = - 2 \tilde u \tilde u_z \tilde v_z +\nu \tilde u_z (\tilde u_z)_{zz},
\label{eqn-uz1}\\
&& \frac{1}{2}({\tilde v}_z^2)_t + \tilde \psi (\tilde v_z^2 )_z
 =  2 \tilde u \tilde u_z \tilde v_z
+\nu \tilde v_z  (\tilde v_z)_{zz}.
\label{eqn-vz1}
\end{eqnarray}
Now, we add (\ref{eqn-uz1}) to (\ref{eqn-vz1}). Surprisingly, 
the nonlinear vortex stretching-like terms cancel each other. We get
\begin{equation}
\left ({\tilde u}_z^2+{\tilde v}_z^2 \right )_t +
2 \tilde \psi \left (\tilde u_z^2 + {\tilde v}_z^2 \right )_z
 = 2 \nu \left ( \tilde u_z (\tilde u_z)_{zz} + \tilde v_z  (\tilde v_z)_{zz}
\right ).
\label{eqn-uvz}
\end{equation}
Further, we note that
\begin{eqnarray*}
\left ( \tilde u_z^2  + \tilde v_z^2 \right )_{zz} =
\left ( 2 \tilde u_z \tilde u_{zz} + 2 \tilde v_z \tilde v_{zz} \right )_z 
= 2 \left ( \tilde u_z (\tilde u_z)_{zz} + \tilde v_z  (\tilde v_z)_{zz}
\right ) + 2 \left [ (\tilde u_{zz})^2 + (\tilde v_{zz})^2 \right ].
\end{eqnarray*}
Therefore, equation (\ref{eqn-uvz}) can be rewritten as
\begin{equation}
\left ({\tilde u}_z^2+{\tilde v}_z^2 \right )_t +
2 \tilde \psi \left (\tilde u_z^2 + {\tilde v}_z^2 \right )_z
 = \nu \left ( \tilde u_z^2  + \tilde v_z^2 \right )_{zz}
- 2 \nu \left [ (\tilde u_{zz})^2 + (\tilde v_{zz})^2 \right ] .
\label{eqn-uvz1}
\end{equation}

Thus, the nonlinear quantity,
$( \tilde u_z^2  + \tilde v_z^2 )$, satisfies a maximum principle
which holds for both $\nu =0$ and $\nu >0$:
\begin{equation}
\|  \tilde u_z^2  + \tilde v_z^2 \|_{L^\infty} \leq
\|  (\tilde u_0)_z^2  + (\tilde v_0)_z^2 \|_{L^\infty} .
\label{maximum_prin}
\end{equation}
Since $\tilde v$ has zero mean, the Poincar\'e inequality
implies that 
$\|\tilde v \|_{L^\infty} \leq C_0 $, with $C_0$ defined by
\[
C_0 = \| \left ( (\tilde u_0)_z^2  + (\tilde v_0)_z^2 
\right )^{\frac{1}{2}}\|_{L^\infty}.
\]
The boundedness of $\tilde u$ follows from the bound on $\tilde v$:
$\|\tilde u (t) \|_{L^\infty} \leq \|\tilde u_0 \|_{L^\infty} \exp ( 2C_0 t)$. 
The higher order regularity follows from the standard estimates.
This proves Theorem 4.

\section{Construction of a family of globally smooth solutions}

In this section, we will use the solution from the 1D model to
construct a family of globally smooth solutions for the 3D
axisymmetric Navier-Stokes equations with smooth initial
data of finite energy. We remark that a special feature of
this family of globally smooth solutions is that the solution
can potentially develop very large dynamic growth and it violates
the so-called smallness condition required by classical global
existence results \cite{CF88,Temam01}.

Let $\bar{u}_1 (z,t)$, $\bar{\omega}_1 (z,t)$, and $\bar{\psi}$
be the solution of the 1D model problem. We will construct a family
of globally smooth solutions of the 3D Navier-Stokes equations
from the solution of the 1D model problem. Denote by
$\tilde{u}(r,z,t)$, $\tilde{\omega}(r,z,t)$ and
$\tilde{\psi}(r,z,t)$ the solution
of the corresponding 3D Navier-Stokes equations. Further, we
define
\begin{equation}
\tilde{u}_1 = \tilde{u}/r, \quad
\tilde{\omega}_1 = \tilde{\omega}/r, \quad
\tilde{\psi}_1 = \tilde{\psi}/r .
\label{eqn-uwp1}
\end{equation}
Let $\phi (r) = \phi_0 (r/R_0)$ be a smooth cut-off
function, where $\phi_0 (r) $ satisfies $\phi_0 (r) = 1 $ if
$0 \leq r \leq 1/2$ and and $\phi_0 (r) = 0 $ if 
$r \geq 1$. Our
strategy is to construct a family of globally smooth functions
$u_1$, $\omega_1$ and $\psi_1$, which are periodic in $z$, such that
\begin{eqnarray}
\tilde{u} &=& r \left (\bar{u}_1 (z,t) \phi (r) + u_1 (r,z,t) \right )
= \bar{u} + u,
\label{utilde}\\
\tilde{\omega} &= &r \left (\bar{\omega}_1 (z,t) \phi (r) +
\omega_1 (r,z,t) \right ) = \bar{\psi} + \psi,\label{wtilde} \\
\tilde{\psi} &= &r \left (\bar{\psi}_1 (z,t) \phi (r) + \psi_1 (r,z,t)
\right ) = \bar{\omega} + \omega ,
\label{psitilde}
\end{eqnarray}
is a solution of the 3D Navier-Stokes equations. 

With the above definition, we can deduce the other two velocity
components
$\tilde{v}^r $ and $\tilde{v}^3$ as follows:
\begin{eqnarray}
&&\tilde{v}^r = -{\tilde \psi}_z = - r \phi (r) \bar{\psi}_{1z}  + v^r (r,z,t)=\bar{v}^r+ v^r , \\
&&\tilde{v}^z = \frac{( r \tilde \psi )_r}{r} = 
\phi (2 \bar{\psi}_1)+ r \phi_r \bar{\psi}_1  + v^z
(r,z,t)=\bar{v}^z+ v^z.
\end{eqnarray}
With the above notations, we can write the velocity vector into two pars as
$\tilde{\textbf{u}}=\bar{\textbf{u}} +\textbf{u}$.

We will choose the initial data for the 1D model of the following form:
\begin{equation}
\bar{\psi}_1 (z,0) = \frac{A}{M^2} \Psi_1 (zM), \;\;
\bar{u}_1 (z, 0 ) = \frac{A}{M} U_1(zM), \;\;
\bar{\omega}_1 (z, 0 ) = A W_1(zM),
\end{equation}
where $A$ and $M$ are some positive constants,
$\Psi_1 (y) $, and $U_1 (y)$ are smooth periodic functions in
$y$ with period 1. Moreover, we assume that $\Psi_1 $
and $U_1$ are odd functions in $y$. Clearly we have
$W_1 = - (\Psi_1)_{yy}$, which is also a smooth,
periodic, and odd function in $y$. It is easy to see that
this feature of the initial data is preserved by the solution
dynamically. In particular,
$\bar{\psi}_1 (z,t)$, $\bar{u}_1 (z,t)$,
$\bar{\omega}_1 (z, t )$ are periodic functions in $z$
with period $1/M$ and odd in $z$ within
each period. Using this property and the
{\it a priori} estimate (\ref{maximum_prin}), we obtain
the following estimate for the solution of the 1D model:
\begin{eqnarray}
&&\|\bar{\psi}_1 (t)\|_{L^\infty} \leq C_0 \frac{A}{M^2},
\label{1D-bound-psi}\\
&&\|\bar{u}_1 (t) \|_{L^\infty} \leq C_0 \frac{A}{M}, \quad
\|\bar{\psi}_{1z} (t) \|_{L^\infty} \leq C_0 \frac{A}{M},
\label{1D-bound-uv} \\ \nonumber \\
&& \|\bar{\omega}_1 (t) \|_{L^\infty} \leq C_0 A ,
\quad
\|\bar{u}_{1z} (t) \|_{L^\infty} \leq C_0 A,
\label{1D-bound-w}
\end{eqnarray}
where
\begin{equation}
C_0 = \| \left ( U_{1y}^2 + W_1^2 \right )^{\frac{1}{2}} \|_{L^\infty}.
\label{amplif-const}
\end{equation}

\vspace{0.1in}
\noindent
{\bf Remark 3.}
As we know from the discussions in the previous sections and
as indicated by (\ref{1D-bound-psi})-\ref{amplif-const}), if the
regularity of the periodic profiles in the initial condition, i.e.,
$U_1$ and $\Psi_1$, is very poor, the solution $\bar{u}_1 (z,t)$ and
$\bar{\psi}_z (z,t)$ will grow very fast dynamically. The amplification
factor is determined by $C_0$ defined in (\ref{amplif-const}).

Let $R_0 = M^{\frac{1}{4}}$. 
From (\ref{1D-bound-psi})-(\ref{amplif-const}) and the definition of
$\bar{\textbf{u}}$, we have
\begin{equation}
\| \bar{\bf u} \|_{L^2} \approx  A R_0^2/M = A/\sqrt{M} ,\quad \| \nabla
\bar{\bf u} \|_{L^2} \approx A R_0^2 = A \sqrt{M} . \label{energy_1-Dpart}
\end{equation}

We would like to emphasize that the corresponding 3D solution
defined by (\ref{utilde}) -(\ref{psitilde}) in general does not
preserve the same special structure in the $z$ direction of the 1D
model problem since the correction terms, $u_1$, $\omega_1$ and
$\psi_1$, are periodic in $z$ with period $1$ instead of period
$1/M$.

We assume that the initial conditions for $u_1$, $\omega_1$, and
$\psi_1$ are chosen in such a way that the principal contributions
to the energy and the enstrophy come from $\bar{\bf u} $,
the mollified solution of
the 1D model. Specifically, we assume that the initial 
condition for ${\tilde {\bf u}}$ satisfies:
\begin{equation}
\| {\tilde {\bf u}}_0 \|_{L^2} \approx A R_0^2/M = A/\sqrt{M} ,\quad
\| \nabla {\tilde {\bf u}}_0 \|_{L^2} \approx A R_0^2 = A \sqrt{M} .
\label{energy_initial}
\end{equation}
Thus, we have
\begin{equation}
\| {\tilde {\bf u}}_0 \|_{L^2} \| \nabla {\tilde {\bf u}}_0 \|_{L^2} \approx A^2 .
\label{energy_prod}
\end{equation}
By choosing $A$ large enough, the above product can be made
arbitrarily large. Thus it violates the classical ``smallness''
condition that guarantees the global existence of the 3D Navier-Stokes
equations \cite{Temam01}.

Furthermore, we have from the energy inequality that
\begin{equation}
\| {\tilde {\bf u}} (t) \|_{L^2} \leq \| {\tilde {\bf u}}_0 \|_{L^2}
\leq A /\sqrt{M} .
\label{energy_cons0}
\end{equation}
Using the above bound and (\ref{energy_1-Dpart}), we obtain {\it a
priori} bound for the perturbed velocity field, ${\bf u}$ in $L^2$
norm:
\begin{equation}
\| {\bf u} (t) \|_{L^2} \leq \frac{A}{\sqrt{M}}.
\label{energy_cons1}
\end{equation}

Let $f = u_1^2$, and define
\begin{eqnarray}
H^2 (t)& = & \int \left ( f^2 +\omega_1^2 \right ) r dr dz =
\int \left ( u_1^4 +\omega_1^2 \right ) r dr dz , \\
E^2 (t) & = & \int \left ( | \nabla f |^2 + | \nabla \omega_1 |^2 \right ) r dr dz ,
\end{eqnarray}
where the integration is over $[0,1]\times [0,\infty)$.

If we further assume that the initial conditions for
$u_1$, $\omega_1$, and $\psi_1$ are odd functions of $z$, then it is
easy to verify that ${\tilde u}_1$, ${\tilde \omega}_1$ and
${\tilde \psi}_1$ are odd functions of $z$ for all times.
Since $\bar{u}_1$, $\bar{\omega}_1$, and $\bar{\psi}_1$ are also
odd functions of $z$, we conclude that
$u_1$, $\omega_1$, and $\psi_1$ are odd functions of $z$ for
all times. It follows by the Poincare inequality that we have
\begin{eqnarray}
&&\int f^2 r dr dz \leq \int f_z^2 r dr dz \leq \int | \nabla f |^2 r dr dz ,\\
&&\int \omega_1^2 r dr dz \leq \int \omega_{1z}^2 r dr dz \leq \int | \nabla
\omega_1 |^2 r dr dz .
\end{eqnarray}
This implies that
\begin{equation}
H \leq E .
\label{eqn-HE}
\end{equation}

Now we can state the main theorem of this section.

\vspace{0.2in}
\noindent
{\bf Theorem 5.} Assume that the initial conditions for
$u_1, \;\omega_1$ and $\psi_1$ are smooth functions of
compact support and odd in $z$. For any
given $A>1$, $C_0>1$ and $\nu>0$,
there exists $C(A, C_0, \nu ) > 0$ such that if
$M > C(A,C_0,\nu )$ and $H(0) \leq 1$, then
the solution of the 3D Navier-Stokes equations given by
(\ref{utilde})-(\ref{psitilde}) remains smooth for all times.

\vspace{0.2in}
\noindent
{\bf Proof.}
First of all, we can use (\ref{eqn-u})-(\ref{eqn-psi}) to derive the
corresponding evolution equations for
$\tilde{u}_1$, $\tilde{\omega}_1$ and $ \tilde{\psi}_1$ as follows:
\begin{eqnarray}
&&(\tilde{u}_1)_t + \tilde{v}^r (\tilde{u}_1)_r +
\tilde{v}^z (\tilde{u}_1)_z = 2 (\tilde{\psi}_1)_z \tilde{u}_1
+ \nu \left ( \tilde{u}_{1zz} + \tilde{u}_{1rr} +
\frac{3 \tilde{u}_{1r}}{r} \right ) ,\label{eqn-u1t} \\
&&(\tilde{\omega}_1)_t + \tilde{v}^r (\tilde{\omega}_1)_r +
\tilde{v}^z (\tilde{\omega}_1)_z =  (\tilde{u}_1^2)_z
+ \nu \left ( \tilde{\omega}_{1zz} + \tilde{\omega}_{1rr} +
\frac{3 \tilde{\omega}_{1r}}{r} \right ) , \label{eqn_w1t}\\
&& - \left ( \tilde{\psi}_{1zz} + \tilde{\psi}_{1rr} +
\frac{3 \tilde{\psi}_{1r}}{r} \right ) = \tilde{\omega}_1 .
\label{eqn-pt1}
\end{eqnarray}

Substituting (\ref{utilde}) into (\ref{eqn-u1t}) and using (\ref{eqn-uwp1}),
we obtain an evolution equation for $u_1$.
\begin{eqnarray}
\frac{\partial u_1}{\partial t} + {\tilde v^r} u_{1r} +
{\tilde v^z} u_{1z}& = & \nu \Delta u_1 + 2{\tilde \psi}_{1z} {\tilde u}_1
-\bar{u}_{1t} \phi - {\tilde v}^r \bar{u}_1 \phi_r \label{u-a} \\
&- & \phi \tilde v^z \bar{u}_{1z}
+ \nu \Delta (\bar{u}_1 \phi ),
\nonumber
\end{eqnarray}
where we have used $\Delta$ to denote the modified Laplacian operator defined by
\[
\Delta w = w_{zz} + w_{rr} + \frac{3 w_{r}}{r} \equiv w_{zz}
+\Delta_{r} w.
\]
On the other hand, we know that $\bar{u}_1$ satisfies the 1D model
equation:
\begin{equation}
\bar{u}_{1t} + 2 \bar{\psi}_1 \bar{u}_{1z}  = \nu 
\bar{u}_{1zz} + 2\bar{\psi}_{1z} \bar{u}_1. \label{u-b}
\end{equation}
Multiplying (\ref{u-b}) by $\phi$ and subtracting the resulting equation
from (\ref{u-a}), we have
\begin{eqnarray}
u_{1t} + {\tilde v^r} u_{1r} + {\tilde v^z} u_{1z}
&= & \nu \Delta u_1 + 2 \left ( {\tilde \psi}_{1z}{\tilde u}_1
- \phi \bar{\psi}_{1z} \bar{u}_1 \right ) \label{eqn-up1} \\
& - & {\tilde v}^r \bar{u}_1 \phi_r - \phi \left ([r \phi_r +2(\phi
-1)] \bar{\psi}_1 + v^z \right ) \bar{u}_{1z} + \nu {\bar u}_1 \Delta_r
\phi . \nonumber
\end{eqnarray}

Similarly, we obtain
\begin{eqnarray}
\omega_{1t} + {\tilde v^r} \omega_{1r} + {\tilde v^z} \omega_{1z}
&= & \nu \Delta \omega_1 + \left ( (u_1 + \bar{u} \phi )_z^2 - \bar{u}_{1z}^2 \phi \right )
\label{eqn-wt1} \\ \nonumber \\
&- & {\tilde v^r} \bar{\omega}_1 \phi_r - \phi \left ([r \phi_r +2
(\phi -1)] \bar{\psi}_1 + v^z \right ) \bar{\omega}_{1z} + \nu
\bar{\omega}_1 \Delta_r \phi . \nonumber
\end{eqnarray}

We divide the analysis into two parts. The first part is devoted to
estimates of the velocity equation. The second part is devoted to
estimates of the vorticity equation.

\vspace{0.1in}
\noindent
{\bf Part I. Estimates for the velocity equation}.
\vspace{0.1in}

First we will present our analysis for the velocity equation. 

Multiply (\ref{eqn-up1}) by $u_1^3$ and integrate over 
$[0,1]\times [0,\infty )$. Using
the incompressibility condition
\[
(r {\tilde v}^r)_r + (r \tilde v^z)_z = 0,
\]
we get
\begin{eqnarray}
\frac{1}{4} \frac{d}{dt} \int u_1^4 r dr dz &\leq &- \frac{3 \nu}{4}
\int \left | \nabla (u_1^2) \right |^2 r dr dz \nonumber \\
&+ & 2 \int u_1^3   \left ( {\tilde \psi}_{1z}{\tilde u}_1
- \phi \bar{\psi}_{1z} \bar{u}_1 \right )
r dr dz
 -  \int {\tilde v}^r \bar{u}_1 \phi_r u_1^3 r dr dz \nonumber \\
&- & \int \phi \left ( [r \phi_r +2\phi (\phi -1)] \bar{\psi}_1 + v^z
\right ) \bar{u}_{1z} u_1^3 r dr dz +
\nu \int \bar{u}_1 ( \Delta_r \phi ) u_1^3 r dr dz \nonumber \\
has the same order
& \equiv & - \frac{3 \nu}{4}
\int \left | \nabla (u_1^2) \right |^2 r dr dz + I + II + III + IV,
\label{est-u}
\end{eqnarray}
where we have used the fact that
\begin{eqnarray}
\int u_1^3 \Delta u_1 r dr dz & = & \int u_1^3 \left (
u_{1zz} + \frac{(r u_{1r})_r}{r} +
\frac{2 u_{1r}}{r} \right ) r dr dz\nonumber \\
& = & -\frac{3}{2} \int \left (u_1^2 u_{1z}^2 + u_1^2 u_{1r}^2 \right )
r dr dz + 2 \int u_1^3 u_{1r} dr dz \nonumber \\
& = & -\frac{3}{4} \int \left [\left ( (u_1^2)_z \right )^2
+ \left ( (u_1^2)_r \right )^2 \right ] r dr dz
- \frac{1}{2} \int u_1^4 (0,z,t) dz \nonumber \\
& \leq & -\frac{3}{4} \int \left [\left ( (u_1^2)_z \right )^2
+ \left ( (u_1^2)_r \right )^2 \right ] r dr dz \;.
\end{eqnarray}

In the following, we will estimate the right hand side of (\ref{est-u}) term by term.

\vspace{0.2in} 
\noindent 
{\bf Estimate for the I-term.} 

Using (\ref{utilde})-(\ref{psitilde}), we have
\begin{eqnarray}
I & = &
2\int u_1^3 \left ( \psi_{1z} u_1 +
\phi \bar{\psi}_{1z} u_1 + \phi \bar{u}_1 \psi_{1z} +
(\phi^2 - \phi ) \bar{\psi}_{1z} \bar{u}_1 \right ) r dr dz \nonumber \\
& \equiv & I_a + I_b + I_c + I_d.
\end{eqnarray}

Using the H\"older inequality, we have
\[
I_a \leq 2 \|\psi_{1z} \|_{L^2} \| f \|_{L^4}^2.
\]
Note that
\begin{eqnarray}
\| \psi_{1z} \|_{L^2}^2 &= &\int \psi_{1z}^2 r dr dz = \int \psi_{1z}^2 d (r^2/2) dz
= - \frac{1}{2} \int r \psi_{1z} \psi_{1zr} r dr dz \nonumber \\
&\leq & \frac{1}{2} \| r \psi_{1z} \|_{L^2} \| \psi_{1zr} \|_{L^2}
\leq \frac{A}{M^{\frac{1}{2}}} \| \psi_{1zr} \|_{L^2} \label{L2-IBP},
\end{eqnarray}
where we have used $r \psi_{1z} =  \psi_{z}$ and (\ref{energy_cons1}) to
obtain
\[
\| r \psi_{1z} \|_{L^2} = \| \psi_{z} \|_{L^2} \leq
\| {\bf u}\|_{L^2} \leq \frac{A}{M^{\frac{1}{2}}}.
\]

On the other hand, using the Sobolev interpolation inequality, we have
\[
\| f \|_{L^4} \leq \| f \|_{L^2}^{\frac{1}{4}} \| \nabla f \|_{L^2}^{\frac{3}{4}} .
\]
This implies that
\[
I_a \leq 2
\| \psi_{1z} \|_{L^2} \|f\|_{L^4}^2 \leq 2 \frac{A^{\frac{1}{2}}}{M^{\frac{1}{4}}}
H E^{\frac{3}{2}} + 
2\frac{(c_2C_0)^{\frac{1}{2}}A}{M^{\frac{11}{8}}} H^{\frac{1}{2}}E^{\frac{3}{2}} ,
\]
where $c_2 = \| \Delta_r \phi_0 \|_{L^\infty}$, and we have used 
\begin{equation}
\| \psi_{1zr} \|_{L^2} \leq \| \omega_1 \|_{L^2}
+ \frac{c_2 C_0 A}{M^{\frac{9}{4}}},
\label{eqn:psi-vort1}
\end{equation}
which we prove in Appendix B.

The estimate for $I_b$ follows from (\ref{1D-bound-uv}):
\[
I_b \leq 2 C_0 \frac{A}{M} \int u_1^4 r dr dz \leq 2 C_0
\frac{A}{M} H^2.
\]
As for $I_c$, we use (\ref{1D-bound-uv}), (\ref{L2-IBP}), and
the H\"older inequality to obtain
\begin{eqnarray}
I_c &\leq & 2 C_0 \frac{A}{M}
\| \psi_{1z} \|_{L^2} \|f \|_{L^3}^{\frac{3}{2}}\nonumber\\
& \leq & 2 C_0 \frac{A}{M} \frac{(A)^{\frac{1}{2}}}{M^{\frac{1}{4}}}
\| \psi_{1zr} \|_{L^2}^{\frac{1}{2}}
 \| f\|_{L^2}^{\frac{3}{4}}
 \| \nabla f\|_{L^2}^{\frac{3}{4}} \nonumber \\
&\leq &  \frac{2C_0 (A)^{\frac{3}{2}}}{M^{\frac{5}{4}}} 
H^{\frac{5}{4}} E^{\frac{3}{4}}
+ \frac{2\sqrt{c_2}C_0^{\frac{3}{2}}A^2}{M^{\frac{19}{8}}}
H^{\frac{3}{4}} E^{\frac{3}{4}},
\end{eqnarray}
where we have used (\ref{eqn:psi-vort1}) and the Sobolev interpolation inequality 
\begin{equation}
\|f\|_{L^3} \leq c_0 \| f\|_{L^2}^{\frac{1}{2}} 
\| \nabla f \|_{L^2}^{\frac{1}{2}}.
\label{Sobolev_L3}
\end{equation}

Finally, we use (\ref{1D-bound-uv}) and the H\"older inequality that
\begin{eqnarray}
I_d \leq 2 C_0^2\frac{A^2}{M^2} \int_{r \leq R_0} |u_1 |^3 r dr dz \leq
2 C_0^2\frac{A^2}{M^2} \left ( \int u_1^4 r dr dz \right )^{\frac{3}{4}}
R_0^{\frac{2}{4}}
 \leq \frac{2C_0^2 A^2}{M^{2-1/8}} H^{\frac{3}{2}}.
\end{eqnarray}
Therefore, we obtain
\begin{eqnarray}
I &\leq & 2 \frac{A^{\frac{1}{2}}}{M^{\frac{1}{4}}}
H E^{\frac{3}{2}} 
+2\frac{(c_2C_0)^{\frac{1}{2}}A}{M^{\frac{11}{8}}} H^{\frac{1}{2}}E^{\frac{3}{2}}
+ 2 C_0 \frac{A}{M} H^2 \nonumber \\
&+& 
\frac{2C_0(A)^{\frac{3}{2}}}{M^{\frac{5}{4}}} H^{\frac{5}{4}} E^{\frac{3}{4}}
+ \frac{2\sqrt{c_2}C_0^{\frac{3}{2}}A^2}{M^{\frac{19}{8}}}
H^{\frac{3}{4}} E^{\frac{3}{4}}
+ \frac{2C_0^2 A^2}{M^{2-1/8}} H^{\frac{3}{2}} .
\end{eqnarray}

\vspace{0.2in}
\noindent
{\bf Estimate for the II-term.}
Using (\ref{1D-bound-uv}), (\ref{energy_cons1}) and the H\"older inequality, we have
\begin{eqnarray}
II &\leq & \frac{c_1 C_0 A}{M R_0} \int | {\tilde v}^r | | u_1^3 | r dr dz
\leq  \frac{c_1 C_0 A}{M R_0}
\| {\tilde v}^r \|_{L^2} \|f\|_{L^3}^{\frac{3}{2}} \nonumber\\
& \leq & \frac{c_1 C_0 A^2}{M^{\frac{3}{2}} R_0} \| f\|_{L^2}^{\frac{3}{4}}
 \| \nabla f\|_{L^2}^{\frac{3}{4}}
\leq  \frac{c_1 C_0A^2}{M^{\frac{7}{4}}} H^{\frac{3}{4}} E^{\frac{3}{4}},
\end{eqnarray}
where $c_1 = \| (\phi_0)_r \|_{L^\infty}$, and we have used the Sobolev 
interpolation inequality (\ref{Sobolev_L3}).

\vspace{0.2in}
\noindent
{\bf Estimate for the III-term.}
Using (\ref{1D-bound-psi}), (\ref{1D-bound-w}), and following the
same steps as in our estimate for the $I_d$-term and the $II$-term, we get
\begin{eqnarray}
III &\leq &(2+c_1) \frac{C_0^2A^2}{M^2} \; \int_{r \leq R_0} | u_1^3 | r dr dz
+ C_0 A  \int | v^z | | u_1^3 | r dr dz\nonumber \\
&\leq & (2+c_1) \frac{C_0^2 A^2}{M^{2-1/8}} H^{\frac{3}{2}}
+  C_0 A \| v^z \|_{L^2} \|f\|_{L^3}^{\frac{3}{2}} \nonumber\\
& \leq &
(2+c_1) \frac{C_0^2 A^2}{M^{2-1/8}} H^{\frac{3}{2}}
+ \frac{C_0 A^2}{M^{\frac{1}{2}}} H^{\frac{3}{4}} E^{\frac{3}{4}}.
\end{eqnarray}

\vspace{0.2in}
\noindent
{\bf Estimate for the IV-term.}

Using (\ref{1D-bound-uv}) and the H\"older inequality, we have
\begin{eqnarray}
IV \leq \frac{\nu c_2 C_0 A}{M R_0^2} \int_{r \leq R_0} | u_1 |^3 r dr dz \leq
\frac{\nu c_2 C_0 A}{M^{\frac{3}{2}}} H^{\frac{3}{2}} R_0^{\frac{2}{4}} \leq
\frac{\nu c_2 C_0 A}{M^{\frac{3}{2}-1/8}} H^{\frac{3}{2}}.
\end{eqnarray}

\vspace{0.1in}
\noindent
{\bf Part II. Estimates for the vorticity equation}.
\vspace{0.1in}

Next, we will present our analysis for the vorticity equation.
Multiplying (\ref{eqn-wt1}) by $\omega_1$ and integrating
over $[0,1]\times[0,\infty)$, we get

\begin{eqnarray}
\frac{1}{2} \frac{d}{dt} \int \omega_1^2 r dr dz & \leq &-
\nu \int | \nabla \omega_1 |^2 r dr dz + \int \left ( {\tilde u}_1^2 -
\bar{u}_1^2 \phi \right )_z \omega_1 r dr dz - \int
{\tilde v}^r \bar{\omega}_1 \omega_1 \phi_r r dr dz  \nonumber\\
&-& \int \phi \left ( [r \phi_r + 2(\phi -1)]\bar{\psi}_1 + v^z
\right ) \bar{\omega}_{1z} \omega_1 r dr dz
+ \nu \int \bar{\omega}_{1} \Delta_r \phi \omega_1 r dr dz \nonumber \\
& \equiv & - \nu \int | \nabla \omega_1 |^2 r dr dz + \bar{I} +
\overline{II} + \overline{III} + \overline{IV} ,
\label{est-w}
\end{eqnarray}
where $ \Delta_r \phi = \phi_{rr} + \frac{3 \phi_r}{r} $.

We will estimate the terms $\bar{I}$ to $\overline{IV}$ one by one.

\vspace{0.2in}
\newpage
\noindent
{\bf Estimate for the $\bar{I}$-term.}

Using (\ref{utilde}) and (\ref{1D-bound-uv}) and integration
by parts, we have
\begin{eqnarray}
\bar{I} &= & -\int (\tilde{u}_1^2 - \bar{u}_1^2 \phi ) \omega_{1z} r dr dz
= - \int (u_1^2 + 2 \bar{u}_1 \phi u_1 +(\phi^2-\phi) \bar{u}_1^2)
\omega_{1z} r dr dz \nonumber \\
&\leq & \int u_1^2 |\omega_{1z}| r dr dz +
\frac{2 C_0A}{M} \int_{r \leq R_0} |u_1| |\omega_{1z}| r dr dz +
\frac{C_0^2 A^2 R_0}{M^2} \| \omega_{1z}\|_{L^2}\nonumber \\
&\leq & \left ( \left (\int u_1^4 r dr dz \right )^{\frac{1}{2}}
+ \frac{2 C_0 A}{M} \left (\int_{r \leq R_0} u_1^2 r dr dz
\right )^{\frac{1}{2}}
+ \frac{C_0^2 A^2}{M^{\frac{7}{4}}} \right ) \| \omega_{1z}\|_{L^2}.
\end{eqnarray}

Let $\Gamma = r {\tilde u}$. It is easy to show that $\Gamma$
satisfies the following evolution equation (see also \cite{BM02}) 
\[
\Gamma_t + {\tilde v}^r \Gamma_r +
\tilde v^z \Gamma_z = \nu (
\Gamma_{zz} + \Gamma_{rr} -
\frac{\Gamma_r}{r} ).
\]
Moreover, for $\tilde u$ smooth, we have $\Gamma |_{r=0} = 0$. 
Thus, $\Gamma$ has a maximum principle, i.e.
\[
\| \Gamma \|_{L^\infty} \leq \| \Gamma_0 \|_{L^\infty} \leq c_0.
\]
This implies that
\[
| r^2 u_1 | \leq | r {\tilde u}|  + r^2 |\phi \bar{u}_1 |
\leq c_0 + R_0^2 \frac{C_0}{M} \leq c_0 +  \frac{C_0}{M^{\frac{1}{2}}}
\leq {\tilde c}_0 .
\]
Therefore, we obtain
\begin{eqnarray*}
\int u_1^4 r dr dz &= & \int (u_1^2 )^2 d (r^2/2) dz
= - \int  r u_1^2 (u_1^2)_r  r d r dz \\
&\leq & \left ( \int r^2 u_1^4 r dr dz \right )^{\frac{1}{2}}
\|\nabla f \|_{L^2}\\
& \leq&
\| r u_1 \|_{L^2}^{\frac{1}{2}}
\left ( \int r^2 u_1^6 r dr dz \right )^{\frac{1}{4}}
\|\nabla f \|_{L^2} \\
&\leq &
{\tilde c}_0^{\frac{1}{4}} \frac{A^{\frac{1}{2}}}{M^{\frac{1}{4}}}
\left ( \int u_1^5 r dr dz \right )^{\frac{1}{4}} \|\nabla f \|_{L^2}.
\end{eqnarray*}
On the other hand, we have
\[
\left ( \int u_1^5 r dr dz \right )^{\frac{1}{4}}
\leq \left ( \int u_1^4 r dr dz \right )^{\frac{1}{8}}
\| f \|_{L^3}^{\frac{3}{8}}  \leq
\left ( \int u_1^4 r dr dz \right )^{\frac{1}{8}}
\| f \|_{L^2}^{\frac{3}{16}} \| \nabla f \|_{L^2}^{\frac{3}{16}}.
\]
Combining the above estimates, we obtain
\begin{eqnarray}
\left ( \int u_1^4 r dr dz \right )^{\frac{1}{2}}
\leq \frac{{\tilde c}_0^{\frac{1}{7}} 
A^{\frac{2}{7}}}{M^{\frac{1}{7}}}
H^{\frac{3}{28}} E^{\frac{19}{28}}.
\label{est-L4}
\end{eqnarray}
Thus, we obtain
\begin{eqnarray}
\bar{I} & \leq & \left (
\frac{{\tilde c}_0^{\frac{1}{7}} 
A^{\frac{2}{7}}}{M^{\frac{1}{7}}}
H^{\frac{3}{28}} E^{\frac{19}{28}}
+ \frac{2 C_0 A }{M} R_0^{\frac{1}{2}} 
\left ( \int u_1^4 r dr dz \right )^{\frac{1}{4}}
+\frac{C_0^2 A^2 }{M^{\frac{7}{4}}} \right ) \| \omega_{1z}\|_{L^2} \nonumber\\
& \leq &
\frac{{\tilde c}_0^{\frac{1}{7}} 
A^{\frac{2}{7}}}{M^{\frac{1}{7}}}
H^{\frac{3}{28}} E^{\frac{47}{28}}+
\frac{2C_0 A }{M^{\frac{7}{8}}} H^{\frac{1}{2}} E
+ \frac{C_0^2 A^2 }{M^{\frac{7}{4}}} E.
\end{eqnarray}

\vspace{0.2in}
\noindent
{\bf Estimate for the $\overline{II}$-term.}

Using  (\ref{1D-bound-w}), (\ref{energy_cons1}), and the H\"older inequality, we get
\begin{eqnarray}
\overline{II} &\leq & \frac{c_1C_0 A}{R_0}
\| {\tilde v}^r \|_{L^2}
 \| \omega_1\|_{L^2} \nonumber \\
&\leq &  \frac{c_1C_0 A}{R_0}
\frac{A}{M^{\frac{1}{2}}}
 \| \omega_1\|_{L^2}  \leq \frac{c_1 C_0 A^2}{M^{\frac{3}{4}}} H .
\end{eqnarray}

\vspace{0.2in}
\noindent
{\bf Estimate for the $\overline{III}$-term.}

Integration by parts gives
\begin{eqnarray}
\overline{III} &= & \int
\phi ([ r \phi_r +2 (\phi -1)]\bar{\psi}_{1z} + v_z^z )
\bar{\omega}_1 \omega_1 r dr dz \nonumber \\
&+ & \int  \phi ( [r \phi_r + 2 (\phi -1)]\bar{\psi}_1 + v^z )
\bar{\omega}_1 \omega_{1z} r dr dz
\end{eqnarray}
We first study the term 
$\int \phi v_z^z \bar{\omega}_1 \omega_{1}  r dr dz$. 
Note that using (\ref{1D-bound-w}), we have
\[
| \int \phi v_z^z \bar{\omega}_1 \omega_{1}  r dr dz |
\leq C_0 A \| (\phi v^z)_z \|_{L^2} \| \omega_1\|_{L^2} .
\]
On the other hand, we have by the Sobolev interpolation
inequality that
\begin{eqnarray*}
\| (\phi v^z)_z \|_{L^2} & \leq & \| \phi v^z \|_{L^2}^{\frac{1}{2}}
\| \nabla (\phi v^z)_z \|_{L^2}^{\frac{1}{2}}\\
&\leq & \frac{\sqrt{A}}{M^{\frac{1}{4}}}
\| \nabla (r \phi v_1^z)_z \|_{L^2}^{\frac{1}{2}}\\
&\leq &  2 \frac{\sqrt{A}}{M^{\frac{1}{4}}} 
 R_0^{\frac{1}{4}}
\| \nabla \omega_1 \|_{L^2}^{\frac{1}{2}}
+  2 \frac{\sqrt{c_2C_0}A }{M^{\frac{7}{8}-1/16}},
\end{eqnarray*}
where we have used
\[
\| \nabla (v_1^z)_z \|_{L^2} \leq  \| \nabla \omega_1 \|_{L^2}
+ \frac{c_2 C_0 A}{M^{\frac{5}{4}}},
\]
which we prove in Appendix B.

Thus we use (\ref{1D-bound-uv}) that
\begin{eqnarray}
\overline{III} &\leq &
(2+c_1)\frac{C_0^2 A^2}{M} R_0 \| \omega_1\|_{L^2}
+ \left (\frac{2C_0 A^{\frac{3}{2}}}{M^{\frac{1}{4}-1/16}} H E^{\frac{1}{2}}
+\frac{2 \sqrt{c_2}C_0^{\frac{3}{2}}A^2}{M^{\frac{7}{8}-1/16}} H \right ) \nonumber\\
&+&(2+c_1) \frac{C_0^2 A^2}{M^2} R_0 \| \omega_{1z}\|_{L^2}
+ C_0 A \| v^z \|_{L^2} \|\omega_{1z}\|_{L^2} \nonumber \\
\nonumber \\
&\leq &
(2+c_1 +2 \sqrt{c_2}) \frac{C_0^2 A^2}{M^{\frac{3}{4}}} H
+\frac{2C_0 A^{\frac{3}{2}}}{M^{\frac{1}{4}-1/16}} H E^{\frac{1}{2}}
+ (3+c_1) \frac{C_0^2 A^2}{M^{\frac{1}{2}}} E .
\end{eqnarray}

\vspace{0.2in}
\noindent
{\bf Estimate for the $\overline{IV}$-term.}

Let
$g(z,t) = \int_0^z \bar{\omega}_1 (\eta, t) d \eta $.
Then we have $g_z = \bar{\omega}_1 $ and $| g | \leq \frac{C_0 A}{M}$.
Thus we have
\begin{eqnarray}
\overline{IV} &= & \nu \int \bar{\omega}_1 (\Delta_r \phi) \omega_1 r dr dz
= \nu \int g_z (\Delta_r \phi )\omega_1 r dr dz \nonumber \\
& = & - \nu \int g (\Delta_r \phi) \omega_{1z} r dr dz  \leq
\frac{\nu c_2 C_0 A}{M R_0^2} \int_{r \leq R_0} | \omega_{1z} | r dr dz \nonumber \\
& \leq & \frac{\nu c_2C_0 A}{M R_0^2} R_0 \| \omega_{1z}  \|_{L^2}
\leq \frac{\nu c_2 C_0 A}{M^{\frac{5}{4}}} E .
\end{eqnarray}

By adding the estimates for $\int u_1^4 r dr dz $ to those for
$\int \omega_1^2 r dr dz $, we obtain an estimate for
$\frac{d}{dt} H^2 $. Note that except for the diffusion
terms, each term in our estimates from $I$ to $\overline{IV}$ can be
bounded  by
\[
\frac{\nu}{16} E^2 + \frac{\epsilon}{16} g (H),
\]
where $g(H)$ is a polynomial of $H$ with positive rational
exponents and positive coefficients
that depend on $C_0$, $A$, and $\nu$, and $\epsilon = \frac{1}{M^\gamma}$
for some $\gamma >0$. Putting all the estimates together,
we get
\begin{eqnarray}
\frac{d}{dt} H^2 \leq  - \frac{\nu}{2} E^2 +
\epsilon g (H) \leq - \frac{\nu}{2} H^2 + \epsilon g (H),
\end{eqnarray}
since $H \leq E$.

For given $A>1$, $C_0 >1$, and $\nu >0$,  we can choose $M$ large enough
so that
\[
- \frac{\nu}{2} + \epsilon g(1) \leq 0 .
\]
Thus, if the initial condition
for $u_1$, $\omega_1$ and $\psi_1$ are chosen such that
$H(0) \leq 1$, then we must have
\[
H(t) \leq 1 , \quad \mbox{for  all } \quad t>0.
\]
Using this ${\it a priori}$ estimate on $H(t)$,
we can easily follow the standard argument to prove the
global regularity of $\psi_1$,  $u_1$ and $\omega_1$ in
higher order norms. This completes the proof of Theorem 5.

\vspace{0.2in}
\centerline{\bf \large Appendix A.}

\vspace{0.2in}
In this appendix, we  prove the following result for the
generalized ODE system.

\vspace{0.2in}
\noindent
{\bf Theorem A.} Assume that ${\tilde u}_0 \neq 0$ and $d \geq 1$. Then the solution
$({\tilde u}(t), {\tilde v}(t))$ of the ODE system (\ref{u-ode1})-(\ref{v-ode1})
exists for all times. Moreover, we have
\begin{equation}
\lim_{t \rightarrow +\infty} {\tilde u} (t) = 0, \quad
\lim_{t \rightarrow +\infty} {\tilde v} (t) = 0 .
\end{equation}
%
%

\vspace{0.2in}
\noindent
{\bf Proof.} We first make a change of variables into the polar coordinate
\footnote{This proof was inspired by a discussion with Mr. Mulin Cheng.}
\begin{equation}
{\tilde v} = r \cos \theta, \quad
{\tilde u} = r \sin \theta.
\end{equation}
Substituting the above change of variables into the ODE system, we obtain
\begin{eqnarray}
r' \cos \theta - r (\sin \theta) \theta' = r^2 \sin^2 \theta - r^2 \cos^2 \theta, \\
r' \sin \theta + r (\cos \theta) \theta' = -d r^2 \cos \theta \sin \theta.
\end{eqnarray}
From the above equations, we can easily derive
\begin{eqnarray}
r' & =& - r^2 \cos \theta \left (\cos^2\theta + (d-1) \sin^2 \theta \right ),
\label{eqn:r}\\
\theta' & = & - r \sin \theta \left ( (d-1) \cos^2 \theta + \sin^2 \theta \right ).
\label{eqn:theta}
\end{eqnarray}
Note that if ${\tilde u}_0 > 0$, then ${\tilde u}(t) > 0$ as long as
$| \int_0^t {\tilde v}(s)ds | < \infty$. Similarly, if
${\tilde u}_0 < 0$, then ${\tilde u}(t) < 0$. Thus, if the solution starts
from the upper (or lower) half plane, it will stay in the upper (or lower)
half plane. Without loss of generality, we may consider the solution
starting from the upper half plane. It follows from (\ref{eqn:theta}) that
$\theta' \leq 0$ since $d \geq 1$, and $r \geq 0$. Therefore, $\theta (t)$ is
monotonically decreasing.  On the other hand, $\theta (t)$ is bounded
from below by zero. As a result, the limit of $\theta (t)$ as
$t \rightarrow \infty$ exists. Let us denote the limiting value
as $\overline{\theta}$. Clearly, we must have
\begin{equation}
\lim_{t \rightarrow \infty} \theta' = 0 , \quad
\lim_{t \rightarrow \infty} \theta = \overline{\theta}.
\label{theta-limit}
\end{equation}

First, we consider the case that the solution starts from the
second quarter (y axis included). We claim that this solution
must cross the $y$-axis into the first quarter. If the solution
stay in the second quarter forever, then $\overline{\theta}$
must be no less than $\pi/2$. From (\ref{eqn:theta}) and
(\ref{theta-limit}), we know that
\begin{equation}
\lim_{t \rightarrow \infty} r = 0.
\label{r-limit}
\end{equation}
However, from (\ref{eqn:r}), we have $r' \geq 0$, which
contradicts with (\ref{r-limit}). The contradiction implies
that the solution must cross the $y$-axis at a later time.

Now we only need to consider the case when the solution
starts from the first quarter since the system is autonomous.
Since $\theta (t)$ decreases monotonically, we obtain
\begin{equation}
r' \leq - \cos^3 (\theta_0 ) r^2 ,
\end{equation}
where we have used the fact that $d \geq 1$ and
$\cos^2 \theta + (d-1) \sin^2 \theta \geq \cos^2 \theta $.
Solving the above ODE inequality gives
\begin{equation}
r (t) \leq \frac{r_0}{1 + r_0 (\cos^3 \theta_0) t} .
\end{equation}

Thus, we conclude that
\begin{equation}
\lim_{t \rightarrow \infty} r (t) = 0 .
\end{equation}

To determine that limiting angle, $\overline{\theta}$, we use
the fact that
\[
\tan (\overline{\theta}) = \lim_{t \rightarrow \infty} \left (
\frac{\tilde u'}{\tilde v'} \right ) =
- \frac{d \tan (\overline{\theta})}{\tan^2 (\overline{\theta}) -1 } .
\]
Since $d \geq 1$, we conclude that $\overline{\theta} = 0$, which
implies
\begin{equation}
\lim_{t \rightarrow \infty} \theta (t) = 0 .
\end{equation}
This completes the proof of Theorem A.

\vspace{0.2in} \centerline{\bf \large Appendix B.}

\vspace{0.2in}

In this appendix, we prove the following two estimates
which relate the $L^2$ norm of the derivatives of $\psi$ to
that of $\omega_1$:
\begin{equation}
\|\psi_{1zz}\|_{L^2} +\|\psi_{1rz}\|_{L^2}+
\|\psi_{1rr}\|_{L^2} + \|\frac{\psi_{1r}}{r} \|_{L^2}
\leq \|w_1\|_{L^2}+
\frac{c_2C_0 A}{M^{\frac{9}{4}}},
\label{B-vort}
\end{equation}
and
\begin{equation}
\|\nabla v_{1z}^z \|_{L^2} = 
\| \nabla ( \frac{2 \psi_{1z}}{r} + \psi_{1rz} )\|_{L^2}
\leq \|\nabla w_1\|_{L^2}+
\frac{c_2C_0 A}{M^{\frac{5}{4}}},
\label{B-dvort}
\end{equation}
where $c_2 = \| \Delta_r \phi_0 \|_{L^\infty}$.

\vspace{0.2in}
\noindent
{\bf Proof.}
From the definition, we have
\[
-\Delta \tilde \psi_1 =\tilde w_1.
\]
Using the definition of $\tilde \psi_1$ and $\tilde w_1$, we 
can rewrite the above equation as
\begin{equation}
-w_1= \Delta \psi_1 + (\Delta_r \phi) \bar{\psi}_1.
 \label{B101}
\end{equation}
Multiplying (\ref{B101}) by $\psi_{1zz}$ and integrating
over $[0,1]\times [0,\infty )$, we obtain:
\begin{eqnarray}
 && \|w_1\|_{L^2}  \|\psi_{1zz}\|_{L^2} \geq \int (\Delta \psi_1
\psi_{1zz} - (\Delta_r \phi) \bar{\psi}_1 \psi_{1zz} )rdrdz \\
&& \geq \int (\psi_{1zz}^2+ \psi_{1rz}^2 ) r dr dz
 - 2\int \psi_{1rz}\psi_{1z} drdz -
\frac{c_2C_0 A}{M^2R_0^2} \int | \psi_{1zz} | r dr dz  
 \\
&& \geq \int (\psi_{1zz}^2+ \psi_{1rz}^2) r dr dz
+ \int_0^1 \psi_{1z}^2(0,z,t)dz -
\frac{c_2C_0 A}{M^2R_0} \| \psi_{1zz}\|_{L^2} \label{estimate-psi1z},
\end{eqnarray}
where we have used (\ref{1D-bound-psi}). This
implies that
\begin{equation}
\|\psi_{1zz}\|_{L^2} +\|\psi_{1rz}\|_{L^2} \leq \|w_1\|_{L^2}
+ \frac{c_2 C_0 A}{M^{\frac{9}{4}}}.
\end{equation}

Next, we multiply (\ref{B101}) by $\Delta_r \psi_{1}$ and integrate
over $[0,1]\times [0,\infty)$. We obtain by using a similar argument
that
\begin{eqnarray}
 && \|w_1\|_{L^2}  \|\Delta_r \psi_{1}\|_{L^2} \geq \int (\Delta \psi_1
\Delta_r \psi_{1} - (\Delta_r \phi) \bar{\psi}_1 \Delta_r \psi_{1} ) rdrdz \\
&& \geq \int [(\Delta_r \psi_{1})^2+ \psi_{1rz}^2]r dr dz
- \frac{c_2 C_0 A}{M^2R_0}\| \Delta_r \psi_{1}\|_{L^2} .
\label{B-2.1}
\end{eqnarray}
On the other hand, we note that
\begin{eqnarray*}
\int (\Delta_r \psi_{1})^2 rdrdz &= &\int ( \psi_{1rr}^2 +9
\frac{\psi_{1r}}{r^2}) rdrdz + 6\int \psi_{1r} \psi_{1rr} dr dz\\
&=&  \int ( \psi_{1rr}^2 +9 \frac{\psi_{1r}}{r^2}) rdrdz
- 3\int_0^1\psi^2_{1r}(0,z,t)dz\\
&=&  \int ( \psi_{1rr}^2 +9 \frac{\psi_{1r}}{r^2}) rdrdz,
\end{eqnarray*}
where we have used the fact that $\psi_{1r}(0,z,t) = 0$ since
$\psi_{1r}$ is odd in $r$. Thus we obtain 
\begin{equation}
\|\psi_{1rr}\|_{L^2} + \|\frac{\psi_{1r}}{r} \|_{L^2}
\leq \|\Delta_r \psi_{1} \| \leq \|w_1\|_{L^2}+
\frac{c_2 C_0 A}{M^{\frac{9}{4}}}.
\label{B-2.2}
\end{equation}
Combining estimate (\ref{B-2.1}) with (\ref{B-2.2}) gives 
the desired estimate (\ref{B-vort}). Similarly, we can 
prove (\ref{B-dvort}).


\vspace{0.2in} 
\noindent 
{\bf Acknowledgments.} We would like to thank Professors
Peter Constantin, Craig Evans, Charles Fefferman, Peter Lax, 
Fanghua Lin, Tai-Ping Liu, Bob Pego, Eitan Tadmor, and S. T. Yau 
for their interests in this work and for some stimulating 
discussions. We also thank Prof. Hector Ceniceros for proofreading 
the original manuscript. The work of Hou was in part supported by NSF 
under the 
NSF FRG grant DMS-0353838 and ITR Grant ACI-0204932 and the work of 
Li was partially supported by NSF grant under DMS-0401174.

\bibliographystyle{amsplain} 
\bibliography{bib}

\end{document}